\documentclass[preprint, 12pt]{elsarticle}
\usepackage{graphicx}
\usepackage{amssymb}
\usepackage{amsmath}
\usepackage{amscd}
\usepackage{mathrsfs}
\usepackage{ntheorem}

\newtheorem{thm}{Theorem}[section]

\newtheorem{lem}[thm]{Lemma}

\newtheorem{defn}[thm]{Definition}

\newtheorem{rem}[thm]{Remark}
\numberwithin{equation}{section}
\newenvironment{proof}[1][Proof]{\textbf{#1.}\ }{\ $\Box$}

\newcommand{\Comp}{\mathbb{C}}

\newcommand{\eps}{\varepsilon}

\newcommand{\A}{\mathcal{A}}

\newcommand{\G}{\mathcal{G}}

\renewcommand{\d}{\mathrm{d}}
\newcommand{\Gal}{\mathrm{Gal}}

\newcommand{\CA}{\G[[\epsilon]]}
\newcommand{\KA}{K\{\A_0(y)\}}

\journal{Annales de l'Institut Henri Poincar\'{e}(C) Analyse Non Lin\'{e}aire}
\begin{document}

\begin{frontmatter}

\title{Nonlinear Differential Galois Theory}

\author{Jinzhi Lei}
\address{Zhou Pei-Yuan Center for Applied Mathematics, 
Tsinghua University,
Beijing, 100084,
P.R.China}

\ead{jzlei@mail.tsinghua.edu.cn}


\begin{abstract}
Differential Galois theory had played important roles in the integrability theory of linear
differential equation. In this paper we will extend the theory to nonlinear differential equation and study the integrability of second order polynomial system. We will propose a definition of the differential Galois group of a differential equation, will study the structure of the group, and will prove the equivalence between the existence of the Liouvillian first integral for the differential equation and the solvability of the corresponding differential Galois group.
\end{abstract}

\end{frontmatter}
\section{Introduction}
\label{sec:0}
In this paper, we will establish the nonlinear differential Galois theory to study the Liouvillian integrability of following nonlinear differential equation
\begin{equation}
\label{eq:1} \dfrac{d x_2}{d x_1} = \dfrac{X_2(x_1,x_2)}{X_1(x_1,x_2)},
\end{equation}
where $X_1$ and $X_2$ are polynomials.  We will propose a definition of the differential Galois group with respect to \eqref{eq:1}, will study the structure of the group, and willl prove the equivalence between the exstience of the Liouillian first integral and the solvability of the group. 

Before we state the main theorem, we give here a brief
outline for the preliminary knowledge of differential
algebra. For detail, refer to \cite{Ka:76} and \cite{Ritt:50}.

\subsection{Preliminary knowledge of differential algebra}
\label{sec:0.4} 

Let $A$ be a ring, by \textit{derivation} of $A$ we mean an
additive mapping $a\mapsto \delta a$ of $A$ into itself satisfying
$$\delta(a b) =(\delta a) b + a (\delta b).$$
We shall say $\delta a$ the \textit{derivative} of $a$. The
\textit{differential ring} $A$ is a commutative ring with unit
together with a derivation $\delta$. If there are $m$ derivations
of $A$, $\delta_i,\ i = 1,2\cdots, m$, satisfying
$$\delta_i \delta_j a = \delta_j \delta_i a,\ \ \forall i,j\in\{1,2,\cdots,m\}, \forall a\in A,$$
we call $A$ together with all the $\delta_i$s a \textit{partial
differential ring}. When $A$ is a field, the \textit{(partial)
differential field} can be defined similarly. In this paper, we
will say differential ring (field) for brevity for both
differential ring (field) and partial differential ring (field).

Let $A$ be any ring, $Y$ be a set of finite or infinite number of
elements. We can generate a ring $A[Y]$ of polynomials of the elements
in $Y$ with coefficients in $A$. In particular, when $A$ is a
differential ring with derivations $\delta_1,\cdots,\delta_m$, and
$Y = \{y_{i_1,i_2,\cdots,i_m}\}\ (i_j = 0,1,\cdots)$ to be the
ordinary indeterminates over $A$, we can extend the derivations of
$A$ to $A[Y]$ uniquely by assigning $y_{i_1\cdots i_{j+1}\cdots
i_m}$ as $\delta_j y_{i_1\dots i_j\cdots i_m}$, and rewriting the
notations as following
$$y_{0\cdots 0} = y,\ \ \ \  y_{i_1\cdots i_m} = \delta_1^{i_1}\cdots \delta_m^{i_m} y.$$
Following the above procedure, we had added a \textit{differential
indeterminate} $y$ to a differential ring $A$. We will denote the
resulting differential ring as $A\{y\}$. The elements of $A\{y\}$ are
\textit{differential polynomials} in $y$. Suppose that $A$ is a
differential field, then $A\{y\}$ is a differential integral
domain, and its derivations can be extended uniquely to the
quotient field. We write $A\langle y\rangle $ for this quotient
field, and its elements are \textit{differential rational
function} of $y$. The notations $\{\ \}$ and $\langle \ \rangle $
will also be used when the adjoined elements are not differential
indeterminates, but rather elements of a larger differential ring
or field.

Let $A$ be any differential ring, then all elements in $A$ with
derivatives $0$ form a subring $C$. This ring is called the ring
of \textit{constants}. If $A$ is a field, so is $C$. Note that $C$
contains the subring that is generated by the unit element of $A$.

Let $A$ be a differential ring, with $\delta_i (i = 1,\cdots,m)$
the derivations. We say an ideal $I$ in $A$ to be a
\textit{differential ideal} if $a\in I$ implies $\delta_i a\in I\
(\forall i)$. An ideal $I$ is said to be a \textit{prime ideal} if
$A\,B\in I$ always implies either $A\in I$ or $B\in I$. Hereinafter,
if not point out particularly, we use the term (prime) ideal in
short for differential (prime) ideal.

Let $A$ and $B$ be two differential rings. A \textit{differential
homomorphism} from $A$ to $B$ is a homomorphism (purely
algebraically) which furthermore commutes with derivatives. The
terms \textit{differential isomorphism} and \textit{differential
automorphism} are self-explanatory.

\subsection{Definitions and statements}
\label{sec:0.2} 
Following first order differential operator $X$ accosiated with the equation \eqref{eq:1} is  convenient in the study,
\begin{equation}
\label{eq:2} X(\omega) = (X_1(x_1,x_2)\delta_1 +
X_2(x_1,x_2)\delta_2)\omega = 0,
\end{equation}
where $\delta_i = \partial/\partial x_i$. From the theory of
differential equation\cite[pp.510-513]{Kaplan:58}, for any
non-critical point $\mathbf{x}^0 = (x_1^0,x_2^0)\in\Comp^2$, the equation \eqref{eq:2} has 
non-constant solution $\omega(x_1,x_2)$
that is analytic at $x^0$. The solution $\omega(x_1,x_2)$ is said
to be a \textit{first integral} of \eqref{eq:1} at $\mathbf{x}^0$.
Furthermore, following lemma can be derived directly from
\cite[Theorem 1, pp. 98]{AB:89}
\begin{lem}
\label{pro:1} Consider the differential equation \eqref{eq:1}, if
$X_1(x_1^0,x_2^0) \not= 0$, and $f(x_2)$ is a function analytic at $x_2 = x_2^0$, then there exists a unique first integral
$\omega(x_1, x_2)$ of \eqref{eq:1}, analytic at $\mathbf{x}^0 =
(x_1^0,x_2^0)$, and
$$\omega(x_1^0, x_2) = f(x_2)$$
for all $x_2$ in a neighborhood of $x_2^0$.
\end{lem}

From the existence of the first integrals of \eqref{eq:1} at
the regular point $\mathbf{x}^0$, we can define the Liouvillian
integrability of \eqref{eq:1} at $\mathbf{x}^0$ as follows.

\begin{defn}
\label{def:3} Let $K$ be the differential field of rational
functions of two variables with derivatives $\delta_1$ and
$\delta_2$, we say $M$ to be a \textit{Liouvillian extension} of
$K$ if there exist $r \geq 0$ and subfields $K_i(i = 0,1,\cdots,
r)$ such that:
$$K = K_0 \subset K_1\subset\cdots\subset K_r = M,$$
where $K_{i+1} = K_i \langle u_i \rangle$, and $u_i\in
K_{i+1}\backslash K_i$ satisfy one of the following:
\begin{enumerate}
  \item $u_i$ is algebraic over $K_i$; or
  \item $\delta_ju_i\in K_i\ \ (j = 1,2)$; or
  \item $\delta_j u_i/u_i\in K_i\ \ (j = 1,2)$.
\end{enumerate}
A function that is contained in some Liouvillian extension of $K$ is
said a \textit{Liouvillian function}.
\end{defn}
\begin{defn}
\label{def:1}Let $K$ be the differential field of rational
functions of two variables, $X$ be defined as \eqref{eq:2}, then
$X$ is \textit{Liouvillian integrable at
$\mathbf{x}^0$} if there exists a first integral $\omega$ of $X$
at $\mathbf{x}^0$, such that $M =K\langle \omega\rangle$ is a
Liouvillian extension of $K$.
\end{defn}

If $X$ is Liouvillian integrable at one point $\mathbf{x}^0\in
\Comp^2$, there exists a first integral that can be obtained from the
rational functions by finite steps of solving algebraic equations, integrals, and exponents of integrals. It is easy to prove by induction that this first integral is analytic in a
dense open set in $\Comp^2$ \cite{S:92}. And hence $X$ it
is also Liouvillian integrable in a dense open set in $\Comp^2$.
Therefore, we can also say that $X$ is Liouvillian integrable.

\begin{defn}
\label{def:2}A group $G$ is \textit{solvable} if there exist a
subgroup series
$$G = G_0\supset G_1\supset\cdots\supset G_m = \{e\}$$
such that for any $0\leq i\leq m-1$, either
\begin{enumerate}
\item $|G_{i}/G_{i+1}|$ is a finite group; or \item $G_{i+1}$ is a
normal subgroup of $G_i$ and $G_i/G_{i+1}$ is an Abelian group.
\end{enumerate}
\end{defn}

Following theorem will be proved in this paper. 

\begin{thm}[Main Theorem]
\label{th:mainthm}
Consider the differential equation \eqref{eq:1}. Assume that $X_1(0,0)\not=0$. 
Let $K$ be the differential field
of rational functions. Then \eqref{eq:1} is Liouvillian integrable if and
only if the differential Galois group of \eqref{eq:1} over $K$ at $(0,0)$
is solvable.
\end{thm}

The keynote in out study is that we don't need to restrict the elements in the differential Galois group to be the automorphoisms. Instead, they are isomorphisms between different extension fields.

\subsection{Historical background}
\label{sec:0.1} 
The first rigorous proof of the non-solvable of a differential equation by quadrature method was given by
Liouville in 1840s\cite{L:41}. Liouville's work was `undoubtedly
inspired by the results of Lagrange, Abel, and Galois on the
non-solvability of algebraic equations by radicals'\cite{Kh:95}.
Since Liouville's pioneer work, many approaches have been
developing toward the integrability theory of differential
equation. The concerning approaches include Lie
group\cite{Olv:99}, monodromy group\cite{Kh:95,Yu:88}, holonomy
group\cite{Ca:00,Ca:01,Ce:98}, differential Galois
group\cite{Ka:76,K:48-2,S:03}, Galois groupoid \cite{Ma:01,
Ma:02}, \textit{etc.}. Let us first recall briefly the subject
of differential Galois theory. For extensive survey, refer to
\cite{Lor:01,S:99}.

At first, we recall Liouville's result. 
Consider following second order linear differential equation
\begin{equation}
\label{eq:17} y'' + a(x) y = 0.
\end{equation}
Liouville proved that the `simple' equation \eqref{eq:17} either
has a solution of `simple' type, or cannot be solved by
quadrature\cite{L:39}. An exposition of Liouville's proof was given in \cite[pp.111-123]{Waston:44}. Explicitly, we have the following result when $a(x)$ is a rational function. 

\begin{thm}[Liouville's Theorem]\cite{L:41}If  $a(x)$ is a rational function, the equation \eqref{eq:17} is solvable by quadrature if and only if it has a solution, $u(x)$ say, such that $u'(x)/u(x)$ is an algebraic
function. 
\end{thm}

From the Liouville's Theorem, to determined the integrability of \eqref{eq:17} with $a(x)$ a rational function, we only need to study the algebraic
function solution of the corresponding Riccati equation (by letting $z = -y'/y$) 
$$z' = z^2 + a(x).$$
One can refer \cite{Kova:86,Lei:00,S:01} for the algorithms to find the Liouvillian solution of \eqref{eq:17} with $a(x)$ a rational function.

Liouville's result was obtained by analytic method. Another approach to the problem
of the integrability of homogenous linear ordinary differential
equation, now known as differential Galois theory, or differential algebra, was established
by Picard and Vessiot at the end of the 19'th century, and well developed by Ritt and Kochin in the next 50 years (see \cite{K:73, Ritt:50} and the references therein). The firm footing step throughout this approach was established by Kolchin in 1948\cite{K:48-1, K:48-2}. For a self-contained exposition of Kochin's work, refer the little book by Kaplansky\cite{Ka:76}.

Consider the linear homogeneous differential equation
\begin{equation}
\label{eq:15} L(y) = y^{(n)} + a_1(x)y^{(n-1)} + \cdots + a_{n-1}(x)y' + a_n(x)y = 0
\end{equation}
with coefficients in a differential $K$.  Let $M$ be an extension of $K$ that contants $n$ linearly independent solutions of \eqref{eq:15}, and has the same field of constants as $K$. Kolchin proved the existence and uniqueness of the extension $M/K$ if $K$ is of characteristic 0 and has an algebraically closed field of constants. This extension is named by Kolchin as the \textit{Picard-Vessiot extension} associated to the linear differential equation. As in th eclassical Galois theory, the differential Galois group of \eqref{eq:15} is defined as the subgroup of all differential automorphisms of $M$  that leaves $K$ elementwise fixed. Kolchin established an isomorphosim between the differential Galois group and an algebraic matric group of degree $n$ over the field of constants   (refer \cite{K:48-2} or \cite[Lemma 5.4]{Ka:76}).  Accordingly, Kolchin was able to prove following connection between the existence of Liouvillian extension and the solvability of the differential Galois group (refer \cite{K:48-2} or \cite[Theorem 5.11-5.12]{Ka:76}).

\begin{thm}[Kolchin's Theorem]\cite{Ka:76, K:48-2}
The linear differential equation \eqref{eq:15} is Liouvillian integrable (the general solution can be obtained by a combination of algebraic functions, integrals and exponentials of
integrals) if and only if the identity component of the corresponding differential Galois group, which is a normal subgroup, is solvable. 
\end{thm}

Kolchin's Theorem is similar to the Galois theory for solvability a polynomial equation by radicals. As application, the previous Liouville's Theorem is easy to be concluded\cite[Theorem 6.4]{Ka:76}. In recent decades, many works by Kovacic, Magid, Mitschi, Singer, Ulmer, \textit{et al.} try to address to both direct and inverse problems
of linear differential Galois theory, for example, see
\cite{Kova:86,Ma:94, Mi:96-1,Mi:96-2,S:93, S:99, S:01, S:03}. The differential Galois theory is also applied to study the non-integrability of Hamiltonian systems\cite{Ruiz:99}. For
details following this approach, one can refer to
\cite{An:94,Ka:76,K:73, Ruiz:99, S:99, S:03} and the references therein.

Besides the linear systems, the solvability of first order
nonlinear differential equation is also interested and will be the main topic of
this paper. Consider the equation
\begin{equation}
\label{eq:10} \dfrac{d x_2}{d x_1} = \dfrac{X_2(x_1,x_2)}{X_1(x_1,x_2)},
\end{equation}
where $X_1$ and $X_2$ are polynomials. The most profound result for
the integrability of this system was obtained by
Singer\cite{S:92} in 1992. Singer proved the following result.

\begin{thm}[Singer's Theorem]\cite{S:92}
If \eqref{eq:10} has a local Liouvillian first integral, then there
is a Liouvillian first integral of the form 
\begin{equation}
\label{eq:11} 
\omega(x_1,x_2) = \int_{(x_1^0,x_2^0)}^{(x_1,x_2)}
R\,X_2\,d x_1 - R\,X_1\,d x_2 
\end{equation}
 where
\begin{equation}
\label{eq:20} R = \exp\left[\int_{(x_1^0,x_2^0)}^{(x_1,x_2)} U \d
x_1 + V \d x_2\right], 
\end{equation} 
with $U$ and $V$ are
rational functions in $x_1$ and $x_2$ such that
$$\partial
U/\partial x_2 =
\partial V/\partial x_1.$$
\end{thm}

Christopher proved that the existence of $R$ in the form \eqref{eq:20} is
equivalent to the existence of an integrating factor of form
$$\exp(D/E)\prod C_i^{l_1},$$
where $D, E$ and the $C_i$ are polynomials in $x_1$ and $x_2$\cite{Ch:99}.
Singer's Theorem was obtained through the method of differential algebra.
The same result was obtained independently by Guan and Lei through Liouville's
approach \cite{Guan:02}. Guan and Lei shown that if
\eqref{eq:10} has local Liouvillian first integral
$\omega(x_1,x_2)$ in the form of \eqref{eq:11}, then
$\delta_i^2\omega/\delta_i\omega\ (i=1,2)$ are rational functions in $x_1$
and $x_2$, here $\delta_i = \frac{\partial\ }{\partial x_i}\ (i =
1,2)$.

Monodromy group of linear differential equation is also important for the integrability of differential equation. In 1970s, Khovanskii proved
that a function is representable by quadrature if and only if its
monodromy group is solvable\cite{Kh:95}. As application of this result, Khovanskii claimed that monodromy group of a linear differential equation is essential for
integrating the equation by quadrature\cite[pp. 128-130,
Khovanskiy's Theorem]{Yu:88}. 

\begin{thm}[Khovanskiy's Theorem]\cite{Yu:88}
If the monodromy group of  a fuchsian system has a solvable normal divisor of finite index, then the system is integrable in quadratures. If the monodromy group does not have this property, then the system is not even integrable by ``generalized quadrature''. This means that the solution of the system cannot be expressed in terms of the coefficients by solving algebraic equations, integration, and composition with entire functions of any number of variables.
\end{thm}

Monodromy group closely relates with the Galois group. The monodromy group of a linear differential equations with regular singular points is Zariski dense in the associated Galois Group\cite{Sch:1887}. \'{Z}o\l\c{a}dek extended the conception of monodromy group  to study the functions defined on $\Comp P^n$ with algebraic
singular set, and to investigate the structure of the monodromy group of the first integrals of a Liouvillian integrable Pfaff equation\cite{Z:98-1}.  Through these studies, \'{Z}o\l\c{a}dek was able to extend Singer's Theorem partly to the integrable
polynomial Pfaff equation\cite{Z:98-1}.

\begin{thm}[Multi-dimensional Singer's theorem] \cite{Z:98-1}
If an integrable polynomial Pfaff equation has a generalized Liouvillian first integral, then it has an integrating factor whose logarithmic differential is a closed rational 1-form.
\end{thm}

In 1990s, the geometry approaches were introduced by
Camacho, Sc\'{a}rdua, \textit{et. al.} to study the Liouvillian
integrability of a nonlinear differential equation. The geometry approaches focus on the 
characters of the foliation associated with the equation. In a series of their works,
the holonomy group that is induced by an invariant algebraic curve of
a polynomial system was studied\cite{Ca:00,Ca:01,Ce:98}. 
Camacho and Sc\'{a}rdua studied the structure of the foliation of a polynomial system with Liouvian first integral. Following result shows such foliations must have simple forms.

\begin{thm}[Camacho-Sc\'{a}rdua Theorem]\cite{Ca:01}
 Let $\mathcal{F}$ be a codimension-one holomorphic foliation on $\mathbb{C}P^n$ admitting a Liouvillian first integral. Assume that one of the algebraic leaves of $\mathcal{F}$ has only non-dicritical infinitesimally hyperbolic singularities. Then either $\mathcal{F}$ is a Darboux foliation or an exponent two Bernoulli foliation after some rational pull-back.
\end{thm}

In 2001, Malgrange published two papers to introduce the
Galois groupoid associated with a foliation with meromorphic
singularities \cite{Ma:01,Ma:02}. For linear differential
equations, Malgrange has showed that this groupoid coincides with the
Galois group of the Picard-Vessiot theory and has proved the required
results in the linear case\cite{Ma:01}. But further development of the
groupoid theory is need to established the theory for nonlinear differential equation\cite{Lor:01}.

The rest of the paper is organized as follows. Section \ref{sec:1} will give the definition of the
differential Galois group, with the discussion of the structure in Section \ref{sec:3}, and leave the proof of Theorem \ref{th:8} to Section \ref{sec:app}. Section \ref{sec:4} will prove the above main theorem. As
application, the differential Galois group of general Riccati and
van der Pol equations will also be discussed in Section \ref{sec:4}. Throughout this paper, 
the base field $K$ always means the field of all
rational functions in $x_1$ and $x_2$ with $\Comp$ the constant field.

\section{Differential Galois group}
\label{sec:1} 
This section will give the definition of the 
differential Galois group of equation \eqref{eq:1} at a
regular point $\mathbf{x}^0 = (x_1^0, x_2^0)$. To this end, we
will first define the group $\G[[\epsilon]]$ that acts at all
first integrals of $X$ at $\mathbf{x}^0$, next studied the
admissible differential isomorphism of \eqref{eq:1} at $x^0$ that
is an element of $\G[[\epsilon]]$, and finally prove that all
admissible differential isomorphisms form a subgroup of
$\G[[\epsilon]]$, which is defined as the
differential Galois group of \eqref{eq:1} at $\mathbf{x}^0$.

Hereinafter, we will assume $\mathbf{x}^0 = (0,0)$ for short. When
we mention a first integral, we will always mean a first integral
that is analytic at $(0,0)$. Following notations will be used
hereinafter. 

Let $\A_0$ denote the set of all functions $f(z)$ of
one variable that are analytic at $z = 0$, and
\begin{eqnarray*}
\A_0^0 &=& \{f(z)\in \A_0|\ \ f(0) = 0\},\\
\A_0^1 &=& \{f(z)\in \A_0^0|\ \ f'(0)\not=0\}.
\end{eqnarray*}
Let $\Omega_{(0,0)}(X)$ denote the set of all first integrals of
\eqref{eq:1} that are analytic at $(0,0)$, and
\begin{eqnarray*}
\Omega_{(0,0)}^0(X) &=& \{\omega(x_1,x_2)\in \Omega_{(0,0)}|\ \
\omega(0,0) = 0\},\\
\Omega_{(0,0)}^1(X) &=& \{\omega(x_1,x_2)\in\Omega_{(0,0)}^0|\ \
\delta_2\omega(0,0)\not=0\}.
\end{eqnarray*}
It is easy to have $f(z) := \omega(0,z)\in \A_0^1$ for any $\omega \in
\Omega_{(0,0)}^1(X)$.

\subsection{The space of all first integrals at a regular point}
\label{sec:1.1} 
Since we will focus on the Liouvillian
integrability of the polynomial system \eqref{eq:1}, it is enough
to concentrate on the first integrals in $\Omega_{(0,0)}^1(X)$ according to the following lemma.
\begin{lem}
\label{pro:2} 
If $X$ is Liouvillian integrable at $(0,0)$, then
there exists a first integral $\omega\in \Omega_{(0,0)}^1$ such
that $M = K\langle \omega\rangle$ is a Liouvillian extension of
$K$.
\end{lem}
\begin{proof}
Let $u$ be a first integral such that $M\langle u\rangle$ is a
Liouvillian extension of $K$. If $u\in \Omega_{(0,0)}^1(X)$, then
the Lemma has been concluded. If
$u\not\in\Omega_{(0,0)}^1(X)$, we can always assume that
$u\in\Omega_{(0,0)}^0(X)$ by subtracting $u(0,0)$ from $u$. Let
$$u(0,x_2) = \sum_{i\geq k} a_i x_2^i \ \ \ \ \ (k\geq 2, a_k\not=0)$$
and
$$f(x_2) = \sqrt[k]{u(0,x_2)} = \left(\sum_{i\geq k} a_i x_2^i\right)^{1/k} = x_2\,\left(\sum_{i\geq k}a_i x_2^{i-k}\right)^{1/k}.$$
Then $f(x_2)\in \A_0^1$. From Lemma \ref{pro:1}, there is a
unique first integral $\omega\in \Omega_{(0,0)}$, such that
$\omega(0, x_2) = f(x_2)$ and therewith $\omega\in
\Omega_{(0,0)}^1(X)$. Moreover, $f(x_2)^k = u(0,
x_2)$ implies $\omega^k = u$, and hence $K\langle \omega\rangle$ is a
Liouvillian extension of $K$.
\end{proof}

From Lemma \ref{pro:1}, there is a one-to-one correspondence
between $\Omega_{(0,0)}^1(X)$ and $\A_0^1$. Hence, it is enough to study the structure of $\A_0^1$. It
is obvious that $\A_0^1$ contains the identity function $e(z) =
z$, and for any $f(z), g(z)\in \A_0^1$, $f\circ g(z) := f(g(z))\in \A_0^1$
and $f^{-1}(z)\in \A_0^1$. Hence $(\A_0^1,\circ)$ is a group. Furthermore, we have the
following result.
\begin{lem}
\label{pro:4} For any $\omega\in\Omega_{(0,0)}^1(X)$, let $\A_0^1(\omega)=\{f\omega)| f\in \A_0^1\}$, then
$$\Omega_{(0,0)}^1(X) =\A_0^1(\omega).$$
\end{lem}
\begin{proof}
It is easy to see that for any $f\in \A_0^1$, $f(\omega)\in \Omega_{(0,0)}^0(X)$. Moreover,
$$\frac{\partial f(\omega)}{\partial x_2}(0,0) = f'(0)\,\delta_2\omega(0,0)\not=0.$$
Hence $f(\omega)\in \Omega_{(0,0)}^1(X)$ and therefore
$\A_0^1(\omega)\subseteq \Omega_{(0,0)}^1(X)$.

For any $\omega, u\in \Omega_{(0,0)}^1(X)$, let
$$g(x_2)= \omega(0,x_2),\ \ h(x_2) = u(0,x_2),$$
then $g,h\in \A_0^1$. Hence $f = h\circ g^{-1}\in \A_0^1$ and
$f(g(x_2)) = h(x_2)$, i.e. $f(\omega(0,x_2)) = u(0,x_2)$. Lemma \ref{pro:1} yields $u = f(\omega)$ and hence
$\Omega_{0,0}^1(X)\subseteq \A_0^1(\omega)$, the Lemma is
proved.
\end{proof}

All elements in $\A_0^1$ map $0$ to $0$. To take the functions that
map $0$ to a nonzero value into account, we adjoin an
infinitesimal variable $\eps$ to the constant field $\Comp$ and
consider the ring of infinite series of $\epsilon$ with
coefficients in $\Comp$. Denote this extension constant ring as
$\Comp[[\epsilon]]$. Then we have the following:
\begin{enumerate}
\item A series $\sum_{i\geq 0}c_i\epsilon^i \in \Comp[[\epsilon]]$ equals $0$ if and
only if all coefficients $c_i$ are $0$; 
\item The derivations $\delta_i$ can be extended to $\Comp[[\epsilon]]$ by setting $\delta_1 \epsilon = \delta_2\epsilon = 0$.
\end{enumerate}

Consider the infinite power series
$$f(z;\epsilon) = \sum_{i,j\geq 0}f_{i,j}z^i\epsilon^j \in \Comp[[z,\epsilon]].$$
The series $f(z;\epsilon)$ is analytic if it is convergent for any
$(z,\epsilon)$ in a neighborhood of $(0,0)$. We can
also write an analytic series $f(z,\epsilon)$ in the form of power series of $\epsilon$ as
$$f(z;\epsilon) = \sum_{i=0}^\infty f_i(z)\,\epsilon^i,$$
where $f_i(z)\in \A_0$. We denote all analytic series as $\A_0[[\epsilon]]$. Let
$$\G[[\epsilon]] = \{f(z;\epsilon) | \ f_0(z)\in \A_0^1\},$$
and define the multiplication in $\G[[\epsilon]]$ as:
$$f(z;\epsilon)\cdot g(z;\epsilon) = f(g(z;\epsilon);\epsilon)$$
for any $ f(z;\epsilon), g(z;\epsilon)\in \G[[\epsilon]]$. Then we have
\begin{lem}
\label{le:46} $(\G[[\epsilon]],\cdot)$ is a group.
\end{lem}
\begin{proof} First, we will show that $(\G[[\epsilon]],\cdot)$ is closure, i.e., for any $f(z;\epsilon), g(z;\epsilon)\in \G[[\epsilon]]$, $f(z;\epsilon)\cdot g(z;\epsilon)\in
\G[[\epsilon]]$. Let
$$f(z;\epsilon) = \sum_{i} f_i(z)\,\epsilon^j,\ \ g(z;\epsilon) = \sum_{i}g_{i}(z)\,\epsilon^j.$$
Then $f(z;\epsilon)$ and $g(z;\epsilon)$ are analytic functions of
$(z,\epsilon)$ at $(0,0)$. Since $g(0,0) = g_0(0) = 0$,
$f(g(z;\epsilon);\epsilon)$ is also an analytic function at
$(0,0)$, i.e., $f(g(z;\epsilon);\epsilon)\in \A_0[[\epsilon]]$.
Moreover, $f(g(z;0);0) = f_0(g_0(z))\in \A_0^1$ and therefore
$f(z;\epsilon)\cdot g(z;\epsilon)\in \G[[\epsilon]]$.

It is easy to verify that
$$(f(z;\epsilon)\cdot g(z;\epsilon)) \cdot h(z;\epsilon) = f(z;\epsilon)\cdot (g(z;\epsilon) \cdot h(z;\epsilon)), $$
Thus $(\G[[\epsilon]],\cdot)$ is associativity.

To prove the identity, we embed $\A_0^1$ into $\A_0^1[[\epsilon]]$ by identifying
$f(z)\in \A_0^1$ with $f(z;0) = f(z) +
\sum_{i\geq0}0\cdot\epsilon^i\in \A_0^1[[\epsilon]]$. Then $e(z;0)
= e(z)\in \G[[\epsilon]]$, and for any $f(z;\epsilon)\in
\G[[\epsilon]]$,
$$e(z;0)\cdot f(z;\epsilon) = e(f(z;\epsilon);0) = f(z;\epsilon),$$
and
$$f(z;\epsilon)\cdot e(z;0) = f(e(z;0);\epsilon) = f(z;\epsilon).$$
Thus, $e(z;0)$ is also an identity of $(\G[[\epsilon]], \cdot)$.

Finally, we only need to prove the invertibility. For any $f(z;\epsilon)\in \G[[\epsilon]]$, we have $(\partial
f/\partial z)(0,0)\not=0$. Thus, the equation $u = f(z;\epsilon)$
has a unique solution $z = f^{-1}(u;\epsilon)$ in the neighborhood
of $(0,0)$ such that $u = f(f^{-1}(u;\epsilon))$. Moreover, we have
$z = f(f^{-1}(z;\epsilon);\epsilon)$ and $z =
f^{-1}(u;\epsilon) = f^{-1}(f(z;\epsilon);\epsilon)$. Thus,
$f^{-1}(z;\epsilon)$ is the inverse element of $f(z;\epsilon)$.
Furthermore, $f^{-1}(z;\epsilon)$ is analytic at $(0,0)$ and
$f^{-1}(z;0)\in \A_0^1$. Thus, we conclude that the inverse element
$f^{-1}(z;\epsilon)\in \G[[\epsilon]]$ and the invertibility is concluded.
\end{proof}

For any $\sigma = f(z;\epsilon)\in \G[[\epsilon]]$ and
$\omega\in\Omega_{(0,0)}^1$, we define the action of $\sigma$
at $\omega$ as $\sigma\omega = f(\omega;\epsilon)$. This is well
defined at the neighborhood of $(0,0)$. Taking account that
$$X(f(\omega;\epsilon)) = X(\sum_{i\geq 0} f_i(\omega)\,\epsilon^i) = \sum_{i\geq 0}X(f_i(\omega))\,\epsilon^i = 0,$$
$\sigma$ maps a first integral $\omega$ to another first integral
$f(\omega;\epsilon)$.

Let $h(z;\epsilon)\in \A_0[[\epsilon]]$ and $f(z;\epsilon)\in
\G[[\epsilon]]$, then
$$h(z;\epsilon)\cdot f(z;\epsilon) = h(f(z;\epsilon);\epsilon)\in\A_0[[\epsilon]]$$
is well defined, and $(h(z;\epsilon)\cdot f(z;\epsilon))\omega =
h(f(\omega;\epsilon);\epsilon)$.

\subsection{Admissible differential isomorphism}
\label{sec:1.2}
For any $\omega\in\Omega_{(0,0)}^1(X)$, an extension field $M$ of
$K$ is obtained by adjoining $\omega$ to $K$, and denoted as $M =
K\langle\omega\rangle$. In this paper, if not mentioned
particularly, $\omega$ will always means a determinate first integral.
In this subsection, we will define and study the admissible
differential isomorphism that is an element of the group $\CA$
with additional restrictions. Throughout this paper, we will call
compactly a map $\sigma:\omega\mapsto\sigma(\omega)$ a
differential isomorphism if there exists a differential
isomorphism from $K\langle\omega\rangle$ to
$K\langle\sigma(\omega)\rangle$ that maps $\omega$ to
$\sigma(\omega)$ with elements in $K$ fixed.

\begin{defn}
\label{def:4}Let $M = K\langle \omega \rangle$ with $\omega\in
\Omega_{(0,0)}^1(X)$. An \textit{admissible differential
isomorphism of $M/K$ with respect to $X$ at $(0,0)$} (a.d.i., singular and plural) is a map $\sigma$ that acts on $M$ with the
following properties: 
\begin{enumerate} 
\item $\sigma$ maps $\omega$ to $f(\omega;\epsilon)$ with some $f(z;\epsilon)\in \CA$;
\item $\sigma: \omega\mapsto f(\omega;\epsilon)$ is a differential
isomorphism;
\item for any $h_i(z;\epsilon)\in \A_0[[\epsilon]]$
($0\leq i \leq m < \infty$), $\sigma$ can be extended to a
differential isomorphism of $K\langle \omega,
h_1(\omega;\epsilon), \cdots, h_m(\omega;\epsilon)\rangle$ that
maps $h_i(\omega;\epsilon)$ to $h_i(f(\omega;\epsilon);\epsilon)$, respectively, with $K$ elementwise fixed.
\end{enumerate} 
\end{defn}

It is obvious that the identity element of $\G[[\epsilon]]$ is an a.d.i.. The following two Lemmas show that the set of all a.d.i. is closure under multiplication and inverse operation.

\begin{lem}
\label{le:19}  If $\sigma,\tau$ are a.d.i., then
$\sigma\cdot\tau$ is an a.d.i..
\end{lem}
\begin{proof}
For any $\varsigma = h(z;\epsilon)\in\A_0[[\epsilon]]$, we will
prove that $\sigma\cdot\tau: \omega\mapsto
(\sigma\cdot\tau)\omega$ can be extended to a differential
isomorphism of $K\langle \omega, \varsigma\omega\rangle$ that maps
$\varsigma\omega$ to $(\varsigma\cdot\sigma\cdot\tau)\omega$ with
$K$ elementwise fixed. 

Since $\tau$ is an a.d.i., $\tau$
is well defined in $K\langle \omega, \sigma\omega,
(\varsigma\cdot\sigma)\omega\rangle$. Hence, the restriction of
$\tau$ at $K\langle \sigma\omega,
(\varsigma\cdot\sigma)\omega\rangle$ is a differential isomorphism
that maps $\sigma\omega$ and $(\varsigma\cdot\sigma)\omega$ to
$(\sigma\cdot\tau)\omega$ and
$(\varsigma\cdot\sigma\cdot\tau)\omega$, respectively, with $K$
elementwise fixed. Consider the following
\begin{center}
\setlength{\unitlength}{0.1 cm}
\begin{picture}(80, 30)
\put(0,25){$K\langle\omega, \varsigma\omega\rangle$}
\put(30,27){$\sigma\cdot\tau$} \put(16,26){\vector(1,0){30}}
\put(47,25){$K\langle (\sigma\cdot\tau)\omega,
(\varsigma\cdot\sigma\cdot\tau)\omega\rangle$}
\put(10,23){\vector(1,-1){18}} \put(16,12){$\sigma$}
\put(24,0){$K\langle \sigma\omega,
(\varsigma\cdot\sigma)\omega\rangle$} \put(45,5){\vector(1,1){18}}
\put(55,10){$\tau|_{K\langle \sigma\omega,
(\varsigma\cdot\sigma)\omega\rangle}$}
\end{picture}
\end{center}
where $\sigma$ and $\tau|_{K\langle \sigma\omega,
(\varsigma\cdot\sigma)\omega\rangle}$ are differential
isomorphisms with $K$ elementwise fixed. Thus
$\sigma\cdot\tau$ is also a differential isomorphism with $K$
elementwise fixed. The extension of $\sigma\cdot\tau$ to
$K\langle\omega,h_1\omega,\cdots,h_m\omega\rangle$ can be proved
similarly. Hence $\sigma\cdot\tau$ is also an a.d.i..
\end{proof}

\begin{lem}
\label{le:33} If $\sigma$ is an a.d.i., then the inverse 
$\sigma^{-1}$ is alo an a.d.i..
\end{lem}
\begin{proof}Similar to the proof of the previous Lemma, it is sufficient to show that for
any $\varsigma = h(z;\epsilon)\in \A_0[\epsilon]$, $\sigma^{-1}$
can be extended to a differential isomorphism of $K\langle
\omega,\varsigma\omega\rangle $ that maps $\varsigma$ to
$(\varsigma\cdot\sigma^{-1})\omega$ with $K$ elementwise fixed.

Consider the extension of $\sigma$ to $K\langle \omega,
\sigma^{-1}\omega, (\varsigma\cdot\sigma^{-1})\omega\rangle$ that
maps $\sigma^{-1}\omega$ and $(\varsigma\cdot\sigma^{-1})\omega$
to $(\sigma^{-1}\cdot\sigma)\omega = \omega $ and
$(\varsigma\cdot\sigma^{-1}\cdot\sigma)\omega = \varsigma\omega$
respectively. Hence, the restricted map $\sigma|_{K\langle
\sigma^{-1}\omega,(\varsigma\cdot\sigma^{-1})\omega\rangle}$ is a
differential isomorphism that maps $K\langle \sigma^{-1}\omega,
(\varsigma\cdot\sigma^{-1})\omega\rangle$ to $K\langle \omega,
\varsigma\omega\rangle$ with $K$ elementwise fixed. Let $\tau =
\left(\sigma|_{K\langle \sigma^{-1}\omega,
(\varsigma\cdot\sigma^{-1})\omega\rangle}\right)^{-1}$, then
$\tau: K\langle \omega, \varsigma\omega\rangle\mapsto K\langle
\sigma^{-1}\omega, (\varsigma\cdot\sigma^{-1})\omega\rangle$ is a
differential isomorphism that maps $\omega$ and $\varsigma \omega$
to $\sigma^{-1}\omega$ and $(\varsigma\cdot\sigma^{-1})\omega$,
respectively, with $K$ elementwise fixed. Thus, we have
$\sigma^{-1} = \tau$ is an a.d.i..
\end{proof}

\subsection{Differential Galois group}
\label{sec:1.3} From the previous discussion, the set of all a.d.i. contains identity element, and satisfies the closure and invertibility, and therefore form a subgroup of $\CA$. This subgroup is our desired differential Galois group.
\begin{defn}
\label{def:5}Let $X$ be defined as \eqref{eq:2}, $K$ be the field
of rational functions, $\omega\in \Omega_{(0,0)}^1(X)$ and $M =
K\langle \omega\rangle$. The \textit{differential Galois group of
$M/K$ with respect to $X$ at $(0,0)$}, denoted as
$\Gal(M/K,X)_{(0,0)}$, is defined as the subgroup of $\CA$ with all elements are admissible differential isomorphism of $M/K$ with respect to $X$ at $(0,0)$.
\end{defn}

Following two Lemmas show that the differential Galois group is
determined `uniquely' by the differential operator $X$ (or the
differential equation \eqref{eq:1}).

\begin{lem}
\label{le:1} Let $u\in \Omega_{(0,0)}^1(X)$ and $N = K\langle u
\rangle$, then $$\Gal(N/K, X)_{(0,0)}\cong \Gal(M/K,X)_{(0,0)}.$$
\end{lem}
\begin{proof}
From Lemma \ref{pro:4}, $u\in \Omega_{(0,0)}^1(x) =
\A_0^1(\omega)$. There is a function $h\in \A_0^1$ such that $u =
h(\omega)$. Let $\tau = h(z;0)\in \CA$, then $u = \tau\omega$,
i.e., $\omega = \tau^{-1}u$.

For any $\sigma\in \Gal(M/K,X)_{(0,0)}$, we will show that
$\tau\cdot\sigma\cdot\tau^{-1}\in \Gal(N/K,X)_{(0,0)}$. To this
end, we only need to show that for any $\varsigma\in
\A_0[[\epsilon]]$, $\tau\cdot\sigma\cdot\tau^{-1}$ can be extended
to a differential isomorphism of $K\langle u, \varsigma u \rangle$
that maps $u$ and $\varsigma u$ to
$(\tau\cdot\sigma\cdot\tau^{-1})u$ and
$(\varsigma\cdot\tau\cdot\sigma\cdot\tau^{-1})u$, respectively,
with $K$ elementwise fixed. 

Since $\sigma\in \Gal(M/K, X)_{(0,0)}$
and $\tau, (\varsigma\cdot\tau)\in \A_0[[\epsilon]]$, $\sigma$ can
be extended to $K\langle \omega, \tau\omega,
(\varsigma\cdot\tau)\omega\rangle$ and maps $\tau\omega = u$ and
$(\varsigma\cdot\tau)\omega = \varsigma u$ to
$(\tau\cdot\sigma)\omega = (\tau\cdot\sigma\cdot\tau^{-1})u$ and
$(\varsigma\cdot\tau\cdot\sigma)\omega =
(\varsigma\cdot\tau\cdot\sigma\cdot\tau^{-1})u$, respectively.
Hence,  $\sigma|_{K\langle u, \varsigma u\rangle}$, the
restriction of $\sigma$ to $K\langle u, \varsigma u\rangle$, is a
differential isomorphism that maps $u$ and $ \varsigma u$ to
$(\tau\cdot\sigma\cdot\tau^{-1})u$ and
$(\varsigma\cdot\tau\cdot\sigma\cdot\tau^{-1})u$, respectively,
with $K$ elements fixed. Thus, we have $\tau\cdot\sigma\cdot\tau^{-1}\in \Gal(N/K,X)_{(0,0)}$. In fact, we have further $\tau\cdot\sigma\cdot\tau^{-1} =
\sigma|_{K\langle u\rangle}$.

Similarly, for any $\eta\in \Gal(N/K,X)_{(0,0)}$,
$\tau^{-1}\cdot\eta\cdot\tau\in \Gal(M/K,X)_{(0,0)}$. Thus
$\sigma\mapsto \tau^{-1}\cdot\sigma\cdot\tau$ is an isomorphism
between $\Gal(M/K,X)_{(0,0)}$ and $\Gal(N/K, X)_{(0,0)}$. The
Lemma has been proved.
\end{proof}

Lemma \ref{le:1} shows that the structure of the differential Galois group of
\eqref{eq:1} at $(0,0)$ is independent to the choice of the first
integral $\omega$. For different choices of first integrals, the corresponding Galois groups are different by a diffeomorphism. Following Lemma will
show that with mild restriction on the regular point, the group is also independent to
the choice of the regular points. We will prove latter (Theorem
\ref{th:8}) that these are all possible cases when the group is of
finite order.

\begin{lem}
\label{le:2} Assume that $\omega\in \Omega_{(0,0)}^1(X)$, $M =
K\langle\omega\rangle$, and $G = \Gal(M/K,X)_{(0,0)}$, we have the
following:
\begin{enumerate}
\item[(1).] If $\omega\in K$, then $\sigma \omega = \omega, \forall
\sigma\in G$; 
\item[(2).] If $(\delta_2\omega)^n \in K$ for some $n\in
\mathbb{N}$, then $\sigma \omega = \mu_n\omega + c(\epsilon),
\forall \sigma\in G$, where $\mu_n$ is a $n$-th root of unity;
\item[(3).] If $\delta_2^2\omega/\delta_2\omega\in K$, then
$\sigma\omega = a(\epsilon)\,\omega + c(\epsilon),\ \forall
\sigma\in G$; 
\item[(4).] If $(2\,(\delta_2\omega)\,(\delta_2^3\omega) -
3\,(\delta_2^2\omega)^2)/(\delta_2\omega)^2\in K$, then
$\sigma\omega = \frac{a(\epsilon)\omega}{1 + b(\epsilon)\,\omega}
+ c(\epsilon),\ \forall \sigma\in G$.
\end{enumerate}
Here $a(\epsilon), b(\epsilon), c(\epsilon)\in
\Comp[[\epsilon]]$, and $c(0) = 0$. Moreover,
for any $(x_1^0, x_2^0)$ such that $X_1(x_1^0, x_2^0)\not=0$ and
$\omega$ is analytic at $(x_1^0, x_2^0)$ and
$\delta_2\omega(x_1^0, x_2^0)\not=0$, the first integral $u$ that is
defined as $u = \omega - \omega(x_1^0, x_2^0)$ is contained in $
\Omega_{(x_1^0, x_2^0)}^1(X)$, and the above results are also
valid for all $\sigma\in \Gal(K\langle u\rangle/K,X)_{(x_1^0,
x_2^0)}$.
\end{lem}
\begin{proof}The first part is proved as follows.

First, (1) is obvious.

(2). Let $(\delta_2\omega)^n  = a\in K$, then
$$a - (\delta_2\omega)^n = 0$$
It is easy to see that for any $\sigma = f(z;\epsilon)\in G$,
$$0  = \sigma (a - (\delta_2\omega)^n) = a - (\delta_2(\sigma(\omega)))^n = a - (\delta_2(f(\omega;\epsilon)))^n = a - f'(\omega;\epsilon)^n(\delta_2\omega)^n$$
Hereinafter, $'$ means the derivative with respect to $z$. Thus,
we have $f'(\omega;\epsilon)^n = 1$ for any $\omega$ in the neighborhood of $\omega(0,0) = 0$, and hence
$f(\omega;\epsilon) = \mu_n\omega + c(\epsilon)$, where
$c(\epsilon)\in \Comp[[\epsilon]]$ and
$\mu_n$ is a $n$-th root of unity.

(3). Let $\delta_2^2\omega/\delta_2\omega = a\in K$, then
$$\delta_2^2\omega - a\,\delta_2\omega = 0.$$
For any $\sigma = f(z;\epsilon)\in G$, we have
\begin{eqnarray*}
0 &=& \sigma(\delta_2^2\omega - a\,\delta_2\omega)\\
&=&\delta_2^2(\sigma\omega) - a\,\delta_2(\sigma\omega)\\
&=&\delta_2^2(f(\omega;\epsilon)) - a\,\delta_2(f(\omega;\epsilon))\\
&=&f''(\omega;\epsilon)\,(\delta_2\omega)^2 +
f'(\omega;\epsilon)\,\delta_2^2\omega
- a\,f'(\omega;\epsilon)\,\delta_2\omega\\
&=&f''(\omega;\epsilon)\,(\delta_2\omega)^2.
\end{eqnarray*}
Hence, we have $f''(\omega;\epsilon) = 0$ and therefore
$\sigma\omega = f(\omega;\epsilon) = a(\epsilon)\,\omega +
c(\epsilon)$ for some $a(\epsilon), c(\epsilon)\in
\Comp[[\epsilon]]$.

(4). Let $(2\,(\delta_2\omega)\,(\delta_2^3\omega) -
3\,(\delta_2^2\omega)^2)/(\delta_2\omega)^2 = a\in K$, then
$$2\,(\delta_2\omega)\,(\delta_2^3\omega) -
3\,(\delta_2^2\omega)^2 - a\,(\delta_2\omega)^2 = 0.$$ 
For any
$\sigma = f(z;\epsilon)\in G$, we have
\begin{eqnarray*}
0&=&\sigma\left(2\,(\delta_2\omega)\,(\delta_2^3\omega) -
3\,(\delta_2^2\omega)^2- a\,(\delta_2\omega)^2\right)\\
&=&2\,\delta_2(\sigma\omega)\,(\delta_2^3(\sigma\omega)) -
3\,(\delta_2^2(\sigma\omega))^2 - a\,(\delta_2(\sigma\omega))^2\\
&=&2\,\delta_2(f(\omega;\epsilon))\,(\delta_2^3(f(\omega;\epsilon)))
- 3\,(\delta_2^2(f(\omega;\epsilon)))^2 - a\,(\delta_2(f(\omega;\epsilon)))^2\\
&=&2\,f'(\omega;\epsilon)\,(\delta_2\omega)\,\left(f'''(\omega;\epsilon)\,(\delta_2\omega)^3
+ 3\,f''(\omega;\epsilon)\,(\delta_2\omega)\,(\delta_2^2\omega) +
f'(\omega;\epsilon)\,\delta_2^3\omega\right)\\
&&{} - 3\,\left(f''(\omega;\epsilon)\,(\delta_2\omega)^2 +
f'(\omega;\epsilon)\,(\delta_2^2\omega)\right)^2 -
a\,\left(f'(\omega;\epsilon)\,(\delta_2\omega)\right)^2\\
&=&2\,f'(\omega;\epsilon)\,f''(\omega;\epsilon)\,(\delta_2\omega)^4
+
6\,f'(\omega;\epsilon)\,f''(\omega;\epsilon)\,(\delta_2\omega)^2\,(\delta_2^2\omega)\\
&&{}+ 2\,(f'(\omega;\epsilon))^2\,(\delta_2\omega)\,(\delta_2^3\omega) - 3\,(f''(\omega;\epsilon))^2\,(\delta_2\omega)^4 \\
&&{}-
6\,f'(\omega;\epsilon)\,f''(\omega;\epsilon)\,(\delta_2\omega)\,(\delta_2^2\omega)
- 3\,(f'(\omega;\epsilon))^2\,(\delta_2^2\omega)^2 -
a\,(f'(\omega;\epsilon))^2\,(\delta_2\omega)^2\\
&=&\left(2\,f'(\omega;\epsilon)\,f''(\omega) -
3\,(f''(\omega;\epsilon))^2\right)\,(\delta_2\omega)^4 \\
&&{}+
(f'(\omega;\epsilon))^2\,(2\,(\delta_2\omega)\,(\delta_2^3\omega)-3\,(\delta_2^2\omega)^2-a\,(\delta_2\omega)^2)\\
&=&(2\,f'(\omega;\epsilon)\,f''(\omega) -
3\,(f''(\omega;\epsilon))^2)\,(\delta_2\omega)^4.
\end{eqnarray*}
Hence, $f(\omega, \epsilon)$ satisfies 
\begin{equation}
\label{eq:foe}
2\,f'(\omega;\epsilon)\,f''(\omega;\epsilon) -
3\,(f''(\omega;\epsilon))^2 = 0,\ \ f(0;\epsilon) = 0.
\end{equation}
The general solution of \eqref{eq:foe} is given by
$$f(\omega;\epsilon) = \frac{a(\epsilon)\,\omega}{1 + b(\epsilon)\,\omega} + c(\epsilon),$$
with $a(\epsilon), b(\epsilon), c(\epsilon)\in
\Comp[[\epsilon]]$.

Finally, since $f(z;\epsilon)\in \G[[\epsilon]]$, we have $f(0;0) = 0$, and hence $c(0) = 0$ in all cases.

For the second part, it is obvious that $u\in \Omega_{(x_1^0,
x_2^0)}^1(X)$, and the above discussions are also valid for $u$. The
proof is complete.
\end{proof}

Similar to classical Galois theory, for any differential subfield
$L$ of $M$ containing $K$, let
$$L' = \{\sigma\in \Gal(M/K,X)_{(0,0)}|\ \sigma a = a, \forall a\in L\}$$
to be the subset of $\Gal(M/K,X)_{(0,0)}$ consisting all a.d.i.
leaving $L$ elementwise fixed. For any subgroup $H$ of $G$, let
$$H' = \{a\in M |\ \sigma a = a,\ \ \forall \sigma\in H\}$$
to be the set of all elements in $M$ left fixed by $H$. Following lemma is obvious from the above definitions.
following results.
\begin{lem}
\label{le:6} Let $L, L_1,L_2$ be subfields of $M$ containing $K$,
and $H,H_1,H_2$ be the subgroups of $G$, then
\begin{enumerate}
  \item[(1).] $L'$ is a subgroup of $G$, and $H'$ is a
  subfield of $M$;
  \item[(2).] $L\subseteq L''$, $H\subseteq H''$;
  \item[(3).] $L_1\supseteq L_2\Rightarrow L_1'\subseteq L_2'$;
  \item[(4).] $H_1\supseteq H_2\Rightarrow H_1'\subseteq H_2'$.
\end{enumerate}
\end{lem}

Let $L$ to be a subfield of $M$ that contains $K$. We can also consider $M$ as an extension field
of $L$ by setting $M = K\langle \omega\rangle = L \langle \omega \rangle$, and the differential Galois group of $M/L$ with respect to $X$ at $(0,0)$ can be defined through the same procedure. We denote this
Galois group as $\Gal(M/L,X)_{(0,0)}$.

\begin{lem}
\label{le:3} Let $L$ be the subfield of $M$ containing $K$, then
$$\Gal(M/L, X)_{(0,0)} = L'.$$ 
In particular, $\Gal(M/K,X)_{(0,0)} = K'$.
\end{lem}
\begin{proof}
It is easy to have $\Gal(M/L,X)_{(0,0)}\subseteq L'$. We will only need to show
that $L'\subseteq \Gal(M/L,X)_{(0,0)}$. 

For $\sigma \in L'$ and $\varsigma\in\CA$, since $K\subseteq L \subseteq K\langle
\omega, \varsigma\omega\rangle$ and $K\subseteq L\subseteq
K\langle \sigma\omega, (\varsigma\cdot\sigma)\omega\rangle$, we
have
$$L\langle \omega, \varsigma\omega\rangle = K\langle \omega, \varsigma\omega\rangle,\
L\langle \sigma\omega, (\varsigma\cdot\sigma)\omega\rangle =
K\langle \sigma\omega, (\varsigma\cdot\sigma)\omega\rangle.$$
From definition \ref{def:4}, $\sigma$ is a differential isomorphism that maps $K\langle \omega, \varsigma \omega \rangle$ onto $K\langle \sigma\omega, (\varsigma\cdot\sigma)\omega\rangle$. Hence, $\sigma$ is also a differential isomorphism that maps $L\langle \omega,
\varsigma\omega\rangle$ onto $L\langle \sigma \omega,
(\varsigma\cdot\sigma)\omega\rangle$, with $L$ elementwise fixed.
From which we conclude that $\sigma\in \Gal(M/L,X)_{(0,0)}$ and
therewith $L'\subseteq \Gal(M/L,X)_{(0,0)}$. The Lemma has been
proved.
\end{proof}

\begin{lem}
\label{le:20}\cite[Lemma 3.1]{Ka:76} Let $M =
K\langle\omega\rangle$, $L$ and $N$ be differential subfields of
$M$ containing $K$ with $N\supset L$, $[N:L] = n$. Let $L'$ and
$N'$ be the corresponding subgroups of $\Gal(M/K,X)_{(0,0)}$. Then
the index of $N'$ in $L'$ is at most $n$.
\end{lem}
\begin{lem}
\label{le:28}\cite[Lemma 3.2]{Ka:76} Let $M =
K\langle\omega\rangle, G=\Gal(M/K,X)_{(0,0)}$ and $H$ and $J$ be
subgroups of $G$ with $H\supset J$ and $J$ of index $n$ in $H$.
Let $H'$ and $J'$ be the corresponding intermediate differential
fields. Then $[J':H']\leq n$.
\end{lem}

\section{Structure of Differential Galois group}
\label{sec:3} 
This section will study the structure of
the differential Galois group $\Gal(M/K,X)_{(0,0)}$. First, we
introduce the preliminary concepts in order to 
describe the structure of the differential Galois group.

\subsection{Quasi-differential polynomial} 
\label{sec:2.2}

Let $y$ be an indeterminate over $K$, and denot by $\A_0(y)$ the
ring
$$\A_0(y) = \{f(y)\ |\ f\in \A_0\}.$$
Adjoining $\A_0(y)$ to $K$ results to a ring $K[\A_0(y)]$ with all
elements of form
$$\sum_{i = 1}^n a_i\,f_i(y),$$
where $a_i\in K, f_i\in \A_0$. The ring $K[\A_0(y)]$
can be extended to a differential ring through the derivatives $\delta_1$ and $\delta_2$
by the rules
$$\delta_i f(y) = f'(y) \delta_iy, \ \ \ (i = 1, 2, f\in \A_0),$$
$$\delta_1(\delta_1^k\delta_2^l y) = \delta_1^{k + 1}\delta_2^l
y,\ \  \delta_2(\delta_1^k\delta_2^l y) = \delta_1^k\delta_2^{l +
1} y ,$$ where $f'\in \A_0$ is a derivative of $f$. Hereinafter, we denote this ring as
$K\{\A_0(y)\}$. 

It is easy to know that all elements in $\KA$ are
polynomials in the derivatives $\delta_1^k\delta_2^l y\ \ (k,l\in
\mathbb{N}_0, k+l > 0)$ with coefficients in $K[\A_0(y)]$. The elements in
$\KA$ differ from differential polynomials in the coefficients that contain not only the polynomials of $y$, but also the terms of form $f(y)$ with $f\in \A_0$. We will call such polynomials of
the derivatives with coefficients in $K[\A_0(y)]$
\textit{quasi-differential polynomials} (QDP, singular and
plural). A \textit{proper quasi-differential polynomial}
(PQDP, singular and plural) is a QDP that involves at least
one proper derivative of $y$. A \textit{regular prime ideal} of
$K\{\A_0(y)\}$ is a prime ideal $\Lambda\subset K\{\A_0(y)\}$ that
contains exclusively PQDP. In this study, we will interest at the regular prime ideal
$\Lambda \subset K\{\A_0(y)\}$ (see Theorem \ref{th:3}).  

The terminologies and results for 
differential polynomials are applicable to PQDP. Let us
recall some basic facts of differential polynomials.  For
detail, refer to \cite{Ritt:50}.

\begin{defn}
\label{def:8}Let
$$w_1 = \delta_1^{i_1}\delta_2^{i_2} y,\ \ \ w_2 = \delta_1^{j_1}\delta_2^{j_2}y,$$
be proper derivatives of $y$, $w_2$ is \textit{higher} then $w_1$ if $j_1
> i_1$ or $j_1 = i_1$ and $j_2 > i_2$. A proper derivative of $y$
is always higher then $y$.
\end{defn}
\begin{defn}
\label{def:9}Let $A$ be a QDP, if $A$ involves proper
derivatives of $y$, by the \textit{leader} of $A$, we mean the
highest of those derivatives of $y$ involved in $A$. If
$A$ involves $y$ but no proper derivatives of $y$, then the leader
of $A$ is $y$. Let $A_1$ be a QDP, and $A_2$ be a PQDP, we
say $A_2$ to be of \textit{higher rank} than $A_1$, if either
\begin{enumerate}
    \item $A_2$ has a higher leader than $A_1$; or
    \item $A_1$ and $A_2$ have the same leader (which is a proper derivative of $y$),
    and the degree of $A_2$ in the leader exceeds that of $A_1$.
\end{enumerate}
Two QDP for which no difference in the rank as created above
will be said to be of the same rank.
\end{defn}
\begin{defn}
\label{def:10} Let $A_1$ be a PQDP, $A_2$ is said to be
\textit{reduced with respect to} $A_1$ if $A_2$ contains no proper
derivative of the leader of $A_1$, and $A_2$ is either zero or of
lower degree than $A_1$ in the leader of $A_1$. Consider a
collection of PQDP
\begin{equation}
\label{eq:51}\Sigma = \{A_1, A_2,\cdots,A_r\},
\end{equation}
if a QDP $B$ is reduced with respect to all the $A_i, (i=
1,\cdots, r)$, then $B$ is said to be reduced with respect to
$\Sigma$.
\end{defn}
\begin{defn}
\label{def:11} Let $F$ be a PQDP with leader $p$, the
QDP $\partial F/\partial p$ is
said the \textit{separant} of $F$. The coefficient of the highest
power of $p$ in $F$ is said the \textit{initial} of $F$.
\end{defn}

\begin{lem}
\label{le:27}\cite[pp.6]{Ritt:50} Let $S_i$ and $I_i$ be, respectively, the separant and initial of
$A_i$ in \eqref{eq:51}, and $F$ be a QDP. There
exist nonnegative integers $s_i,t_i, i= 1,\cdots,r$, such that
when a suitable linear combination of the $A$ and their
derivatives is subtracted from
$$S_1^{s_1}\cdots S_r^{s_r}I_1^{t_1}\cdots I_r^{t_r} F,$$
the remainder is reduced with respect to \eqref{eq:51}.
\end{lem}

Let $\Lambda\in K\{A_0(y)\}$ be a regular prime idea, $X(y) =
X_1\delta_1y + X_2\,\delta_2y\in \Lambda$, and $\{X(y)\}$ is the
differential ring that generated by $X(y)$. Let $A(y)\in \Lambda$
with the lowest rank and irreducible. If
$\Lambda\varsupsetneq\{X(y)\}$, it is easy to see that $A(y)$ involves
no $\delta_1y$ and its derivatives. Let $\delta_2^ry$ be the
leader of $A(y)$. Then $A(y)$ is a polynomial of the derivatives
$\delta_2y,\delta_2^2y,\cdots,\delta_2^ry$, with coefficients
$A_i(x_1,x_2,y)\in K[\A_0(y)]$. From Lemma \ref{le:27}, the regular
prime idea $\Lambda$ is the least regular prime idea containing
$X(y)$ and $A(y)$. Thus, according to \cite[pp. 4-5]{Ritt:50}, the \textit{characteristic set} of
$\Lambda$ consists of $A(y)$ and $X(y)$. We will see latter that the number $r$ is
important to determine the structure of $\Lambda$, and named as
the \textit{order} of $\Lambda$, denoted by $\mathrm{ord}(\Lambda) = r$. If
$\Lambda = \{X(y)\}$, then the characteristic set of $\Lambda$
contains only one element $X(y)$, and the order is said to be
$\infty$.

\subsection{Structure of Differential Galois Group}
\label{sec:3:1}

\begin{lem}
\label{le:21} If there exists $A(y)\in K[\A_0(y)]$\
($A(y)\not\equiv 0$), and $\omega\in \Omega_{(0,0)}^1(X)$, such
that $A(\omega(x_1,x_2)) = 0$ for all $(x_1,x_2)$ in a neighborhood of $(0,0)$, then $K$ contains a first integral of $X$.
\end{lem}
\begin{proof}
Hereinafter, we will write $A(\omega) \equiv 0$ in short for $A(\omega(x_1,x_2)) = 0$ for all $(x_1,x_2)$ in a neighborhood of $(0,0)$. 

Let
$$\Sigma_0 = \{A(y)\in K[\A_0(y)]|\ A(\omega) \equiv 0,\ \ A(y)\not\equiv 0\}.$$
Then $\Sigma_0\not=\emptyset$. For any $A(y)\in \Sigma_0$, we can write $A(y)$ in the form as
\begin{equation}
\label{eq:Ay}
A(y) = \sum_{k=1}^n \alpha_k f_k(y),
\end{equation}
with $\alpha_k\in K,\ f_k(y)\in  \A_0(y)$. It is not unique to express $A(y)$ with the form \eqref{eq:Ay}.  Within all possible expressions, there is one with the shortest length $n$. We call this shortest length $n$ the length of $A(y)$, and denote by $n(A)$. 

Let $A(y)$ to be the element in $\Sigma_0$ with the
smallest length. If $n(A) =  1$, then $A(y) = \alpha_1 f_1(y)$, and hence
$f_1(\omega) \equiv 0$, i.e., $\omega$ is a constant. Therefore we have $n
> 1$ since the first integral $\omega$ cannot be a constant.

When $n > 1$, write
$$A(y) = \alpha_1 f_1(y) + \alpha_2 f_2(y) + \cdots + \alpha_n f_n(y)$$
with $n = n(A)$, then
$$A(\omega) = \alpha_1 f_1(\omega) + \alpha_2 f_2(\omega) + \cdots + \alpha_n f_n (\omega) \equiv 0,$$
and hence
$$X (A(\omega)) = X (\alpha_1) f_1(\omega)  + X(\alpha_2) f_2(\omega) + \cdots +
X(\alpha_n) f_n(\omega) \equiv 0.$$ Therefore
\begin{equation}
\label{eq:axa}
\alpha_1 X(A(\omega)) - X(\alpha_1) A(\omega) = \sum_{i = 2}^n
(\alpha_1 X(\alpha_i) - X(\alpha_1) \alpha_i) f_i(\omega) \equiv 0.
\end{equation}
If $\alpha_1 X(\alpha_i) - X(\alpha_1)\alpha_i\not=0$ for some $i$, then \eqref{eq:axa} implies that
$$B(y) = \sum_{i = 2}^n(\alpha_1 X(\alpha_i) - X(\alpha_1) \alpha_i) f_i(y)$$
is an element in $\Sigma_0$ with smaller length then 
$A(y)$, which is contradict. Thus, we have
\begin{equation}
\label{eq:7} \alpha_1 X(\alpha_i) - X(\alpha_1) \alpha_i = 0,\ \ (i=2,3,\cdots, n).
\end{equation}
Moreover, it is easy to see that $\alpha_2/\alpha_1$ is not a constant and $X(\alpha_2/\alpha_1) = 0$ by \eqref{eq:7}. Hence $\alpha_2/\alpha_1$ is a first integral of $X$ contained in $K$, and the Lemma is proved.
\end{proof}

From Lemma \ref{le:21}, if $K$ contains no first integral of $X$,
and there is a $A(y)\in K\{\A_0(y)\}$ such that $A(\omega) \equiv 0$ for some
$\omega\in \Omega_{(0,0)}^1(X)$, then $A(y)\not\in K[\A_0(y)]$, i.e., $A(y)$ is a PQDP that involves some proper derivatives of $y$.

Let $\Lambda$ be e regular prime ideal of PQDP, we say a function $u(x_1,x_2;\epsilon)$ satisfies $\Lambda$ if for any
$F(x_1,x_2,y,\delta_1y, \delta_2y,\cdots)\in \Lambda$, while
substitute $y$ and the derivatives in $F$ with $u$ and the
corresponding derivatives, the resulting expression is zero for
all $x_1, x_2$ and $\epsilon$ that are small enough.

\begin{thm}
\label{th:3} Let $K$ and $X$ be defined as previous. Assume that $K$ contains no
first integral of $X$. Let $M = K\langle \omega\rangle$ with
$\omega\in \Omega_{(0,0)}^1(X)$. Then there exists a regular prime
ideal $\Lambda$ of PQDP such that:
\begin{enumerate}
  \item[(1).] For every $\sigma_f\in \Gal(M/K,X)_{(0,0)}$, let $\sigma_f\omega = f(\omega;\epsilon)\ (f(z;\epsilon)\in \G[[\epsilon]])$, then
$f(\omega(x_1,x_2);\epsilon)$ satisfies $\Lambda$.
  \item[(2).] Given $f(z;\epsilon)\in \CA$ such that $f(\omega(x_1,x_2);\epsilon)$ satisfies $\Lambda$,
  there exists $\sigma_f\in \Gal(M/K,X)_{(0,0)}$ such that $\sigma_f(\omega) = f(\omega;\epsilon)$.
\end{enumerate}
\end{thm}
\begin{proof}
Let $y$ be a differential indeterminate over $K$, define the
natural homomorphism from $K\{\A_0(y)\}$ to $K\{\A_0(\omega)\}$
that maps $h(y)$ to $h(\omega)$\ ($\forall h\in \A_0$). Let
$\Lambda$ to be the kernel of the homomorphism, then $\Lambda$ is
a regular prime ideal of $K\{\A_0(y)\}$. We will prove that
$\Lambda$ fulfil the requirement of the Theorem.

(1). Let $\sigma_f\in \Gal(M/K,X)_{(0,0)}$ and $\sigma_f\omega =
f(\omega;\epsilon)$. Then $f(z;\epsilon)\in \G[[\epsilon]]$. For
any $F(y)\in \Lambda$, i.e., $F(\omega) \equiv 0$, there exist
$h_i\in\A_0\ (i = 1,\cdots, m)$, such that $F(y)\in K\{y,
h_1(y),\cdots,h_m(y)\}$. Write $F(y)$ as a differential
polynomial of $h_1(y), \cdots, h_m(y)$, i.e.,
$$F(y) = F(y, h_1(y),\cdots,h_m(y)),$$
then
$$F(\omega, h_1(\omega),\cdots,h_m(\omega)) \equiv 0.$$
Since $\sigma_f\in \Gal(M/K,X)_{(0,0)}$, $\sigma_f$ can be
extended to a differential isomorphism of $K\{\omega, h_1(\omega),
\cdots, h_m(\omega)\}$ that maps $\omega$ and $h_i(\omega)$ to
$f(\omega;\epsilon)$ and $h_i(f(\omega;\epsilon))$, respectively.
Thus, we have
$$F(f(\omega;\epsilon), h_1(f(\omega;\epsilon)), \cdots,
h_m(f(\omega;\epsilon))) \equiv 0.$$ 
The requirement (1) has been proved.

(2). Now, consider $f(z;\epsilon)$ in $\CA$ such that
$f(\omega;\epsilon)$ satisfies $\Lambda$, we will show that $\sigma_f\in \Gal(M/K,X)$. For any
$h_1(z;\epsilon), \cdots, h_m(z;\epsilon)\in \A_0[[\epsilon]]$,
consider the maps
$$\begin{array}{cccc}
\pi: &K\langle y, h_1(y;\epsilon)), \cdots,
h_m(y;\epsilon)\rangle&\mapsto &
K\langle \omega, h_1(\omega;\epsilon), \cdots, h_m(\omega;\epsilon)\rangle\\
&y&\mapsto&\omega\\
&h(y;\epsilon)&\mapsto& h(\omega;\epsilon)\\
&\delta_jy&\mapsto&\delta_j\omega
\end{array}
$$
and
$$\begin{array}{cccc}
\pi_\sigma: &K\langle y, h_1(y;\epsilon), \cdots,
h_m(y;\epsilon)\rangle&\mapsto &
K\langle f(\omega;\epsilon), h_1(f(\omega;\epsilon);\epsilon), \cdots, h_m(f(\omega;\epsilon);\epsilon)\rangle\\
&y&\mapsto& f(\omega;\epsilon)\\
&h(y;\epsilon)&\mapsto& h(f(\omega;\epsilon);\epsilon)\\
&\delta_jy&\mapsto&\delta_j\omega
\end{array}
$$
where $h\in \A_0(y)$ and $j = 1,2$. Let the kernels of $\pi$ and
$\pi _\sigma$ be $\Gamma$ and $\Gamma_\sigma$, respectively.

We will show that $\Gamma = \Gamma_\sigma$.
Then $K\langle \omega, h_1(\omega;\epsilon), \cdots,
h_m(\omega;\epsilon)\rangle$ is isomorphic to $K\langle
f(\omega;\epsilon), h_1(f(\omega;\epsilon);\epsilon), \cdots,
h_m(f(\omega;\epsilon);\epsilon)\rangle$ with the isomorphism
$\sigma:\omega\mapsto f(\omega;\epsilon),
h_i(\omega;\epsilon)\mapsto h_i(f(\omega;\epsilon);\epsilon)$.
Therefore, $\sigma$ is an admissible differential isomorphism, and the Theorem is proved. 

First, we will prove $\Gamma \subseteq \Gamma_\sigma$. For any $F(y,h_1(y;\epsilon),\cdots,h_m(y;\epsilon))\in \Gamma$, we write $F$ in the form of the power series in $\epsilon$
$$
F(y,h_1(y;\epsilon),\cdots,h_m(y;\epsilon)) = \sum_{i =
0}^{\infty}F_i(y,h_{i,1}(y),\cdots,h_{i,m_i}(y))\,\epsilon^i\ \ \
(h_{i,j}\in \A_0)
$$ where $F_i$ are differential polynomials. Then
$$F(\omega,h_1(\omega;\epsilon),\cdots,h_m(\omega;\epsilon)) = \sum_{i =
0}^{\infty}F_i(\omega,h_{i,1}(\omega),\cdots,h_{i,m_i}(\omega))\,\epsilon^i\equiv 0$$
i.e., $F_i(\omega,h_{i,1}(\omega),\cdots,h_{i,m_i}(\omega)) \equiv 0$ for all $i$.
Thus, the coefficients $F_i$ are contained in $\Lambda$. Now,
assume that $f(\omega;\epsilon)$ satisfies $\Lambda$, then
$$F_i(f(\omega;\epsilon),h_{i,1}(f(\omega;\epsilon)),\cdots,h_{i,m_i}(f(\omega;\epsilon))) \equiv 0,\ \ \ (\forall i).$$
Thus,
$$F(f(\omega;\epsilon),h_1(f(\omega;\epsilon)),\cdots,h_m(f(\omega;\epsilon))) \equiv 0,$$
and hence $F(y,h_1(y;\epsilon),\cdots,h_m(y;\epsilon))\in
\Gamma_\sigma$. Therefore, $\Gamma\subseteq\Gamma_\sigma$.

Next, we will show that $\Gamma_\sigma\subseteq \Gamma$. If on the
contrary, there exists $F(y;\epsilon)\in \Gamma_\sigma$ but
$F(y;\epsilon)\not\in \Gamma$, then
$$F(f(\omega; \epsilon);\epsilon) =
F(f(\omega;\epsilon),h(f(\omega;\epsilon);\epsilon),\cdots,h_m(f(\omega;\epsilon);\epsilon)) \equiv 0,$$ 
but
$$F(\omega;
\epsilon) = F(\omega,
h_1(\omega;\epsilon),\cdots,h_m(\omega;\epsilon)) \not\equiv 0.$$
Write $F(y;\epsilon)$ as a power series in $\epsilon$
$$F(y;\epsilon) = \sum_{i = 0}^{\infty}F_i(y)\,\epsilon^i,\ \ \ (F_i(y)\in K\{\A_0(y)\}),$$
and let $k$ the smallest index such that $F_i(y)\in \Lambda$
for any $0\leq i \leq k-1$ and $F_k(y)\not\in\Lambda$. From the
assumption that $f(\omega;\epsilon)$ satisfies $\Lambda$, we have
$F_i(f(\omega;\epsilon)) = 0$ for any $0 \leq i\leq k-1$. Thus,
let $f_0(z) = f(z;0)$, we have
$$F(f(\omega;\epsilon);\epsilon) = \epsilon^k\,F_k(f(\omega;\epsilon)) +
\sum_{i \geq  k+1}F_i(f(\omega;\epsilon))\,\epsilon^{i} =
F_k(f_0(\omega))\,\epsilon^k + h.o.t. \equiv 0$$
and therefore $F_k(f_0(\omega)) \equiv 0$.

Now, we obtain a $F_k(y)\not\in\Lambda$, and $F_k(f_0(\omega)) \equiv
0$. Let $A(y)$ in $\Lambda$ with the lowest rank, and hence $A(y)$
and $X(y)$ make up the characteristic set of $\Lambda$. Let $S(y)$
and $I(y)$ to be the separant and initial of $A(y)$, respectively (If
$\Lambda = \{X(y)\}$, we take $S(y) = I(y) = X_2(x_1,x_2)$). It is
clear that $S(y), I(y) \not\in \Lambda$. Lemma \ref{le:27} yields that
there exist nonnegative integrals $s,t$, and $R(y)\in
K\{\A_0(y)\}$ that is reduced with respect to $\Lambda$, such that
$$S(y)^s I(y)^t F_k(y) - R(y) \in \Lambda.$$
Since $\Lambda$ is a prime ideal and $S(y), I(y),
F_k(y)\not\in\Lambda$, we have $S(y)^sI(y)^tF_k(y)\not\in
\Lambda$, and thus $R(y)\not=0$. From the above discussion, 
$f_0(\omega)$ satisfies both $\Lambda$ and $F_k(y)$, and hence
$R(f_0(\omega)) = 0$, which implies $R(f_0(y))\in \Lambda$.
Whereas, simple computation shows that $R(f_0(y))$ has the same
rank as $R(y)$, and therefore is reduced with respect to $\Lambda$,
which is contradict. Hence we have proved $\Gamma_\sigma \subseteq
\Gamma$.

Now, we have concluded $\Gamma = \Gamma_\sigma$, and the Theorem has been proved.

\end{proof}

\begin{rem}
\label{rem:3} Theorem \ref{th:3} indicates that when $X$ has no first integral in $K$, the regular prime ideal $\Lambda$ is essential to determine the differential Galois group. Here we associate the order of the differential Galois group with that of $\Lambda$ as following.
\begin{enumerate}
\item Recall the order of $\Lambda$ that was defined in
Section \ref{sec:2.2}. If $A(y)$ and $X(y)$ make up the
characteristic set of $\Lambda$, and the highest derivative of
$A(y)$ is $\delta_2^ry$ ($0<r<+\infty$), then
$\mathrm{ord}(\Lambda) = r$. If $\Lambda = \{X(y)\}$, then
$\mathrm{ord}(\Lambda) = \infty$. It is easy to obtain from Theorem
\ref{th:3} that if $\mathrm{ord}(\Lambda) = \infty$, then
$\Gal(M/K,X)_{(0,0)} = \G[[\epsilon]]$. 
\item If there exists
$g\in \A_0^1$ such that $g(\omega)=u$ is contained in $K$, then $g(y)-
u\in\Lambda$. In this case, we say the order of $\Lambda$ is
$0$. 
\item We will also define the order of the differential
Galois group as the order of $\Lambda$.
\end{enumerate}
\end{rem}

\begin{thm}
\label{th:8} Let $\Lambda$ the
prime ideal in Theorem \ref{th:3} and $r = \mathrm{ord}(\Lambda)$, then either $0\leq r\leq 3$ or
$r = \infty$. Moreover, let $G = \Gal(M/K,X)_{(0,0)}$, we have the following
\begin{enumerate}
\item[(1).] If $r = 0$, then $K$ contains a first integral of $X$, and
\begin{equation}
\label{eq:16}G = \{e\}.
\end{equation} 
\item[(2).] If
$r = 1$, then there exists $\omega\in\Omega_{(0,0)}^1(X)$ such
that $(\delta_2\omega)^n\in K$ for some $n\in \mathbb{N}$. Let $M
= K\langle\omega\rangle$, then
\begin{equation}
\label{eq:12} G = \{f(z;\epsilon)\in
\G[[\epsilon]]\ \Big|\ f(z;\epsilon) = \mu\,z + c(\epsilon),\ \
c(0) = 0, \mu^n = 1\}.
\end{equation} 
\item[(3).] If $r = 2$, then there exist $\omega\in
\Omega_{(0,0)}^1(X)$ such that $\delta_2^2 \omega /\delta_2
\omega\in K$. Let $M = K\langle\omega\rangle$, then
\begin{equation}
\label{eq:13} G =
\{f(z;\epsilon)\in\G[[\epsilon]]\ \Big|\ f(z;\epsilon) =
a(\epsilon)\,z + c(\epsilon),\ \ c(0) = 0\}.
\end{equation}
\item[(4).] If $r = 3$, then there exists $\omega\in
\Omega_{(0,0)}^1(X)$ such that
$$\dfrac{\delta_2 \omega\cdot\delta_3\omega - 3\,\delta_2^2\omega}{(\delta_2 \omega)^2}\in K.$$
Let $M = K\langle\omega\rangle$, than
\begin{equation}
\label{eq:14} G =
\{f(z;\epsilon)\in\G[[\epsilon]]\ \Big|\ f(z;\epsilon) =
\frac{a(\epsilon)\,z}{1 + b(\epsilon)\,z} + c(\epsilon),\ \ c(0) =
0\}.
\end{equation}
\item[(5).] If $r = \infty$, then $G =
\G[[\epsilon]]$.
\end{enumerate}
Here $a(\epsilon), b(\epsilon), c(\epsilon)\in \Comp[[\epsilon]]$. In particular, if $G$ is solvable, then $r\leq 2$.
\end{thm}
We leave the proof of Theorem 3.9 to Section \ref{sec:app}. From Theorem
\ref{th:8} and Lemma \ref{le:2}, the differential Galois group is
independent to the choice of the point $(0,0)$. Moreover, from
Lemma \ref{le:1}, the structure of the group is also independent to the
choice of the first integral $\omega$. Hence, we can omit the
point $(0,0)$ and the particular extension $M$, and simply say $\Gal(M/K,X)$ the differential Galois group of $X$ over $K$. This group is determined uniquely by the equation \eqref{eq:1} or the operator $X$ and will tell the insight of the integrability of the differential equation.

\section{Liouvillian integrability of the polynomial system}
\label{sec:4} 
We are now ready to prove the main theorem of this paper. 

\subsection{Preliminary results of Galois theory}
\label{sec:4.1}
\begin{defn}
\label{def:13}Let $K$ be a (differential) field, $M$ be an
extension field of $K$, $G$ be a set of isomorphisms of $M$, with
$K$ elementwise fixed. $M$ is \textit{normal} over $K$ with
respect to $G$ if there is no element in $M\backslash K$ that is fixed by all elements in $G$.
\end{defn}

Obviously, we have
\begin{lem}
\label{le:18}Let $G=\Gal(M/K,X)_{(0,0)}$ and $H$ be a subgroup of $G$. Let
$$H' = \{a\in M|\ \sigma a = a, \forall \sigma\in H\},$$
then $M$ is normal over $H'$ with respect to $H$.
\end{lem}

\begin{lem}
\label{le:22} If $K$ contains no first integral of $X$, then for
any $\omega\in \Omega_{(0,0)}^1(X)$, $M = K\langle \omega\rangle $
is normal over $K$ with respect to $\Gal(M/K,X)_{(0,0)}$.
\end{lem}
\begin{proof}
We only need to prove that for any $\alpha\in M\backslash K$,
there exists $\sigma\in \Gal(M/K,X)_{(0,0)}$, such that
$\sigma\alpha\not=\alpha$.

Let $\alpha = p(\omega)/q(\omega)$, with $p(y), q(y)\in
K\{\A_0(y)\}$. Without loss of generality, we assume further that
$p(y), q(y)$ are reduced with respect to the prime ideal $\Lambda$ given by Theorem \ref{th:3}.
Therefore, $p(\omega)\not=0$ and $q(\omega)\not=0$. Write $p(y)$ and
$q(y)$ explicitly as
$$p(y) = p(\mathbf{x},y,\delta_2y,\cdots,\delta_2^{r-1}y),\ \ q(y) = q(\mathbf{x},y,\delta_2y,\cdots,\delta_2^{r-1}y)$$
where $\mathbf{x} = (x_1, x_2)$, $r = \mathrm{ord}(\Lambda)$, and let
$$A(\mathbf{x},y) = p(\mathbf{x},y,\delta_2y,\cdots,\delta_2^{r-1}y) -
\alpha(\mathbf{x})\,q(x_1,x_2,y,\delta_2y,\cdots,\delta_2^{r-1}y).
$$

Since $q(\omega)\not\equiv0$, there exists
$\mathbf{x}^0 = (x_1^0, x_2^0)$ such that $\omega(\mathbf{x})$ is
analytic at $\mathbf{x}^0$, and $q(\omega(\mathbf{x}^0))\not=0$.
Let $\mathbf{c} =
(\omega(\mathbf{x}^0),\cdots,\delta_2^{r-1}\omega(\mathbf{x}^0))$
and $\mathbf{x}^*$ in a neighborhood $U$ of
$\mathbf{x}^0$ such that
$$q(\mathbf{x}^*, \mathbf{c})\not=0\ \ \ \mathrm{and}\ \ \ \ A(\mathbf{x}^*, \mathbf{c})\not=0.$$
We claim that such $\mathbf{x}^*$ always exists. If on the
contrary, we have
$$A(\mathbf{x},\mathbf{c}) = p(\mathbf{x},\mathbf{c}) - \alpha(\mathbf{x})\,q(\mathbf{x},\mathbf{c}) = 0$$
for any $\mathbf{x}\in U$ such that
$$q(\mathbf{x}, \mathbf{c})\not=0,$$
then
$$\alpha(\mathbf{x}) = \frac{p(\mathbf{x},\mathbf{c})}{q(\mathbf{x},\mathbf{c})}\in K,$$
which is contradict to the fact that $\alpha\in M\backslash K$.

Let $\mathbf{c^*} =
(\omega(\mathbf{x}^*),\cdots,\delta_2^{r-1}\omega(\mathrm{x}^*))$,
and $\epsilon$ to be an infinitesimal parameter, then 
$$q(\mathbf{x}^*,\mathbf{c^*}+ \epsilon\,(\mathbf{c} - \mathbf{c^*}))\not=0$$
and
$$A(\mathbf{x}^*,\mathbf{c^*}+ \epsilon\,(\mathbf{c} - \mathbf{c^*}))\not=0.$$
And hence
$$\frac{p(x_1^*,x_2^*,\mathbf{c^*}+ \epsilon\,(\mathbf{c} - \mathbf{c^*}))}
{q(x_1^*,x_2^*,\mathbf{c^*}+ \epsilon\,(\mathbf{c} -
\mathbf{c^*}))}\not=\alpha(\mathbf{x}^*).$$ 
By Theorem \ref{th:8},
it is easy to verify that when the order $r\not=0$, there exists
$\sigma\in\Gal(M/K,X)_{(0,0)}$ such that
$$(\delta_i\sigma\omega)(\mathbf{x}^*) = c_i^* + (c_i-c_i^*)\epsilon\ \ (i = 0,1,\cdots,r-1).$$
Therefore
$$(\sigma\alpha)(\mathbf{x}^*)  = \frac{p(\mathbf{x}^*,\mathbf{c^*}+ \epsilon\,(\mathbf{c} - \mathbf{c^*}))}
{q(\mathbf{x}^*,\mathbf{c^*}+ \epsilon\,(\mathbf{c} -
\mathbf{c^*}))}\not= \alpha(\mathbf{x}^*),$$ 
i.e., $\sigma\alpha\not=\alpha$. The Lemma has been proved.
\end{proof}

In following Lemmas, we let $L,N,M$ be the extension fields of $K$, with
$M = K\langle \omega\rangle $, $K\subset L\subset N \subset M$,
and $G = \Gal(M/K,X)_{(0,0)}$. Assume that $N = L\langle u\rangle
$ with $u$ satisfying $\delta_i u\in L$ or $\delta_iu/u\in L,\ (i
=1,2)$. Then $L'$ and $N'$ are subgroups of $G$, and $N''$ is a
subfield of $M$.
\begin{lem}
\label{le:13} Any $\sigma\in L'$ maps $N$ to $N''$.
\end{lem}
\begin{proof}We have $\sigma a = a$ for any $\sigma$ in $L'$ and $a$ in $L$. At first, assume that $N=\langle u\rangle$. If $\delta_iu = a_i\in L, (i = 1,2)$, then $\delta_i (\sigma
u) = a_i,\ (i =1,2)$. Thus $\sigma u = u + c(\epsilon) $ with
$c(\epsilon)\in \Comp[\epsilon]$, and therewith $\sigma u\in N''$,
which implies $\sigma N \subseteq N''$. 

The proof for the
case $\delta_i u/u \in K$ is similar by the fact that $\sigma u =
c(\epsilon)\,u$ with $c(\epsilon)\in \Comp[[\epsilon]]$ for any
$\sigma$ in $L'$.
\end{proof}
\begin{lem}
\label{le:14} $N'$ is a normal subgroup of $L'$, and $L'/N'$ is Abelian.
\end{lem}
\begin{proof}Let $\sigma\in L', \tau\in N'$, then $\sigma a \in
N''$ for any $a\in N$, and therefore $\tau (\sigma a) = \sigma a$.
Thus $(\sigma^{-1}\cdot \tau\cdot \sigma) a = \sigma^{-1} (\sigma
a) = a$, and hence $\sigma^{-1}\cdot \tau\cdot \sigma\in N'$. This
implies that $N'$ is a normal subgroup of $L'$.

Assume that $N = L\langle u\rangle$ with $\delta_iu\in L\
(i=1,2)$.  From the proof of Lemma \ref{le:13}, for any $\sigma\in
L'$, we have $\sigma u = u + c(\epsilon)$ for some $c(\epsilon)\in
\Comp[[\epsilon]]$. Thus, the subgroup $H = \{\sigma|_N\ |\ \sigma\in L'\}$ is
isomorphic to a subgroup of the addition group $\Comp[[\epsilon]]$
and hence is Abelian. Consider the homomorphism from $L'$ to
$H$ that maps $\sigma$ to $\sigma|_N$. The kernel of the map is
$N'$, and the image is $H$. Thus, $L'/N'$ is isomorphic to $H$ and
is Abelian.

The case that $N = L\langle u\rangle$ with $\delta_iu/u\in L\
(i=1,2)$ can be proved similarly.
\end{proof}

From Lemma \ref{le:20} and \ref{le:14}, we have:
\begin{lem}
\label{le:31} Let $L,N,M$ be the extension fields of $K$, with $M =
K\langle \omega\rangle , N = L\langle u\rangle $, and $K\subset
L\subset N \subset M$. Assume
further that $K, L, N$ contain no first integral of $X$. We have
\begin{enumerate}
    \item[(1).] if $u$ is algebraic over $L$, then $|L'/N'|\leq [N:L]$; and
    \item[(2).] if $\delta_i u \in L$ or $\delta_i u/u\in L\ (i = 1,2)$,
    then ${N'}$ is a normal subgroup of ${L'}$, and $L'/N'$ is
    Abelian.
\end{enumerate}
\end{lem}

\begin{lem}
\label{le:35} Assume that $M = K\langle \omega \rangle$ is normal
over $L$ with respect to $G$. If for every $\sigma\in G$, there
exist $c(\epsilon)\in \Comp[[\epsilon]]$, such that
$$\sigma \omega = \omega + c(\epsilon),$$
then $M$ is a Liouvillian extension of $L$.
\end{lem}
\begin{proof}
For any $\sigma\in G$, we have
$$\sigma (\delta_i\omega) = \delta_i\omega,\ \ (i = 1,2).$$
Since $M$ is normal over $L$ with respect to $G$, we have
$\delta_i\omega\in L, (i = 1,2)$, and $M$ is a Liouvillian
extension of $L$.
\end{proof}

\begin{lem}
\label{le:30} Assume that $M = K\langle \omega \rangle$ is normal
over $K$ with respect to $G$. If for every $\sigma\in G$, there
exist $a(\epsilon), c(\epsilon)\in \Comp[[\epsilon]]$, such that
$$\sigma \omega = a(\epsilon) \omega + c(\epsilon),$$
then $M$ is a Liouvillian extension of $K$.
\end{lem}
\begin{proof}
For any $\sigma\in G$, we have
$$\sigma (\delta_i^2\omega/\delta_i\omega) = \delta_i^2\omega/\delta_i\omega,\ \ (i = 1,2).$$
Since $M$ is normal over $K$ with respect to $G$, there exist
$a_i\in K$, such that
$$\delta_i^2\omega = a_i \delta_i\omega, \ \ (i = 1,2).$$
Taking account that $X_1\delta_1\omega + X_2\delta_2\omega = 0$,
there exists $\mu\in M$ such that $\delta_1\omega = \mu X_2,
\delta_2\omega = -\mu X_1$, and hence
$$\dfrac{\delta_1\mu}{\mu} = a_1 - \dfrac{\delta_1 X_2}{X_2} \in K,\ \ \
\dfrac{\delta_2\mu}{\mu} = a_2 - \dfrac{\delta_2 X_1}{X_1}\in K.$$
Thus $M$ is a Liouvillian extension of $K$, and
$$K\subseteq K\langle\mu\rangle\subseteq K\langle \omega \rangle = M.$$
\end{proof}

\subsection{Proof of the Main Theorem}
\label{sec:4.2}

\begin{proof}[Proof of Theorem \ref{th:mainthm}]
(1). Let $G$ to be the differential Galois group of \eqref{eq:1} over $K$ at $(0,0)$. If $K$ contains a first integral of $X$, then $G$ contains
exclusively the identity mapping, and is solvable.

Now, we assume that $K$ contains no first integral of $X$, and $X$
is Liouvillian integrable. Then there exists $\omega\in
\Omega_{(0,0)}^1(X)$ such that $M = K\langle \omega\rangle $ is a
Liouvillian extension of $K$. From definition \ref{def:3},
suppose that:
$$K = K_0\subset K_1\subset K_2\subset \dots\subset K_m = M,$$
with $K_{i+1} = K_i\langle u_i\rangle $, where either $u_i$ is
algebraic over $K_i$ or $\delta_j u_i\in K_i$ or $\delta_j
u_i/u_i\in K_i$ ($j = 1,2$). Let $G_0 = \Gal(M/K,X) ( = K'), G_i =
K_i'$, $(i = 1,2,\cdots,m)$, then
$$G = G_0\supseteq G_1 \supseteq G_2\supseteq\cdots\supseteq G_m = \{e\}.$$
By Lemma \ref{le:31}, either $|G_i/G_{i+1}|\leq[K_{i+1}:K_i] <
\infty$ or $G_{i+1}$ is a normal subgroup of $G_i$, and
$G_i/G_{i+1}$ is Abelian. Thus, $G$ is solvable by definition \ref{def:2}.

 (2). If the differential Galois group of \eqref{eq:1} over $K$ at $(0,0)$ is solvable,
by Theorem \ref{th:8}, either $K$ contains a first integral of $X$,
or the Galois group has order $r = 1$ or $r = 2$.

By Theorem \ref{th:8}, if $r = 1$, then there exists $\omega\in
\Omega_{(0,0)}^1(X)$ and $n\in \mathbb{N}$ such that
$$G = \{f(z;\epsilon)\in \G[[\epsilon]]\ |\ f(z;\epsilon) = \mu\, z + c(\epsilon), \ c(\epsilon)\in \Comp[[\epsilon]], c(0) = 0,\ \mu^n = 1\},$$
where $M = K\langle\omega\rangle$ and $G = \Gal(M/K,X)_{(0,0)}$.
Let 
$$G_0 = \{\sigma\in G|\ \sigma\omega = \omega + c(\epsilon),\  c(\epsilon)\in \Comp[[\epsilon]], c(0) = 0\},$$
then $G_0$ is a subgroup of $G$, and $|G/G_0| = n$. By Lemma
\ref{le:22}, $K = G'$. Hence, according to Lemma \ref{le:28},
$$[G_0':K] = [G_0':G']\leq [G/G_0] = n,$$
which means that $G_0'$ is an algebraic extension of $K$. By Lemma
\ref{le:18}, $M$ is normal over $G_0'$ with respect to $G_0$.
Hence Lemma \ref{le:35} is applicable and yields that $M$ is a
Liouvillian extension of $G_0'$, and consequently a Liouvillian
extension of $K$.

If $r= 2$, then there exists a first integral $\omega\in
\Omega_{(0,0)}^1$ such that
 $$G  = \{f(z;\epsilon)\in \G[[\epsilon]]\ |\ f(z;\epsilon) = a(\epsilon)\, z + c(\epsilon),\ c(0) = 0\},$$
where $M = K\langle \omega \rangle$ and $G = \Gal(M/K,X)_{(0,0)}$ as in previous. Again, by Lemma \ref{le:22}, $K = G'$. Hence Lemma \ref{le:30} is
applicable, and $M$ is a Liouvillian extension of $K$. The Theorem
has been proved.
\end{proof}

\subsection{Applications}

From the proof of Lemmas \ref{lem:a.3}-\ref{lem:a.15} in Section \ref{sec:app}, the explicit
method to determine the differential Galois group can be given as follows.
\begin{thm}
\label{th:2}Consider the differential equation \eqref{eq:1}, let
\begin{equation}
\label{eq:Bi}
B_i = -X_1\delta_2^{i+1}(\frac{X_2}{X_1}),\  \ \ i = 0,1,2
\end{equation}
and $r$ to be the order of the corresponding differential Galois
group, then
\begin{enumerate}
\item[(1).] $r=0$ if and only if $K$ contains a first integral of $X$; 
\item[(2).] $r=1$ if and only if $K$ contains no first integral of $X$, and there exists $a\in K, a\not=0$, and $n\in
\mathbb{N}$, such that
\begin{equation}
\label{eq:31} X(a) = n\,B_0\,a.
\end{equation}
 \item[(3).] $r=2$ is and only if \eqref{eq:31} is not satisfied by all $a\in K$ and $n\in \mathbb{N}$, and there exists
$a\in K$, such that 
\begin{equation} 
\label{eq:32} 
X(a) = B_0\,a + B_1.
\end{equation}
 \item[(4).] $r=3$ if and only if \eqref{eq:32} is not satisfied by all $a\in K$, and there exists $a\in K$, such that
\begin{equation}
\label{eq:33} X(a) = 2\,B_0\,a +  B_2. 
\end{equation}
 \item[(5).] $r=\infty$ if and only if \eqref{eq:33} is not satisfied by all $a\in K$.
\end{enumerate}
\end{thm}

It is easy to see that the equation
\begin{equation}
\frac{d y}{d x} = p(x)
\end{equation} 
with $p(x)$ a polynomial has the order of the differential Galois group $r = 0$. The general homogenous linear equation
\begin{equation}
\frac{d y}{d x} = p(x) y
\end{equation}
has the order $r = 1$, and the general nonhomoegenous linear equation
\begin{equation}
\frac{d y}{d x} = p(x) y + q(x)
\end{equation}
has the order $r = 2$. Following result shows that the general Riccati equation is an example with order $r = 3$.
\begin{thm}
\label{cor:3} The differential Galois group of the general Riccati
equation
\begin{equation}
\label{eq:riccati}
\frac{d y}{d x} = p_2(x)\,y^2 + p_1(x)\,y + p_0(x)
\end{equation}
has order $r = 3$.
\end{thm}
\begin{proof} We have known that the general Riccati equation \eqref{eq:riccati} does not have Liouvillian first intergral\cite{L:41}, and hence the Theorems \ref{th:mainthm} and \ref{th:8} indicate that either the order $r =3$ or $r = \infty$. 

From the equation \eqref{eq:riccati}, we have $X_1 = 1$ and $X_2 = p_2(x)\,y^2 + p_1(x)\,y + p_0(x)$
and $\delta_2 = \partial/\partial y$, and $B_2 = 0$ from \eqref{eq:Bi}. Thus, the equation \eqref{eq:33} has solution $a = 0$, and the order is $3$.
\end{proof}

Next, we will show an example with the order of the differential Galois group to be infinity.

Consider the van der Pol equation
\begin{equation}
\label{vdp} \left\{\begin{array}{rcl} \dot{x}_1 &=& x_2 - \mu
(\dfrac{x_1^3}{3} - x_1),\\
\dot{x}_2 &=& -x_1 \end{array}\right.
 \ \ \ (\mu \not = 0).
\end{equation}
The van der Pol equation is well known for the existence of a limit cycle. 
Following Lemma was proved independently by Cheng et al.\cite{Cheng:95} and Odani\cite{Oda:95}, respectively, at almost the same time. 

\begin{lem} \cite{Cheng:95, Oda:95}
\label{lem:vdp} The system of the van der Pol equation \eqref{vdp} has no algebraic solution curves. In particular, the limit cycle of it is not algebraic.
\end{lem}

\begin{thm}
\label{cor:2} 
The order of the differential Galois group of the van der Pol equation \eqref{vdp} is
infinity.
\end{thm}
\begin{proof}
Let
$$X_1(x_1,x_2) = x_2 - \mu
(\dfrac{x_1^3}{3} - x_1),\ \ X_2(x_1,x_2) = -x_1,$$ 
the equation
\eqref{eq:33} for the van der Pol equation \eqref{vdp} becomes
\begin{equation}
\label{eq:37} X_1^3\,  X(a) + 2 x_1 X_1^2\,  a + 6 x_1 = 0.
\end{equation}
We will only need to prove that \eqref{eq:37} has no rational
function solution. 

If \eqref{eq:37} has a rational function
solution $a = a_1/a_2$, where $a_1, a_2$ are relatively prime
polynomials, then $a_1$ and $a_2$ satisfy
$$X_1^3\,(a_2\, X(a_1) - a_1\, X(a_2)) - 2 x_1\,X_1^2\,a_1\,a_2 + 6\, x_1\, a_2^2 = 0.$$
Hence, there exist a polynomial $c(x_1,x_2)$, such that
\begin{eqnarray}
\label{eq:38} X_1^3\, X(a_2) &=& c\, a_2,\\
\label{eq:39} X_1^3\, X(a_2) &=& (c - 2\, x_1\, X_1^2)\,a_1 - 6
x_1\,a_2.
\end{eqnarray}

Let $a_2 = X_1^k\, b_2$, where $k \geq 0$, $b_2$ is a nonzero
polynomial and $(b_2, X_1) = 1$. Substitute $a_2$ into
\eqref{eq:38} yields
$$X_1^3\, X(b_2) + k\,  X_1^2\,  X(X_1)\,  b_2 = c\,b_2.$$
Thus, $b_2| X_1^3\,  X(b_2)$ and therewith $b_2 | X(b_2)$, i.e.,
either $b_2$ is a constant or $b_2(x_1,x_2)=0$ is an algebraic solution curve of
\eqref{vdp}. However, Lemma \ref{lem:vdp} has shown that \eqref{vdp} 
has no algebraic solution curves. Thus $b_2$ must be a constant. Let $b_2 = 1$ without loss of generality , we have
\begin{equation}
\label{eq:40} a_2 = X_1^k,\ \ c = k\, X_1^2\, X(X_1).
\end{equation}
Substitute \eqref{eq:40} into \eqref{eq:39} yields
\begin{equation}
\label{eq:41} X_1^3\, X(a_1) = (k\,X_1^2\, X(X_1) - 2\, x_1\,
X_1^2) a_1 - 6\, x\, X_1^k,\ \ \ (k\geq 0).
\end{equation}
Note that
$$(k\, X(X_1) - 2 x_1) = -k\, \mu\, (x_1^2 - 1)\, X_1 - (k + 2)\, x_1,$$
\eqref{eq:41} becomes
$$X_1^3\, X(a_1) = X_1^2\,(-k\, \mu\, (x_1^2 - 1)\, X_1 - (k + 2)\, x_1)\, a_1 - 6\, x_1\, X_1^k.$$
If $k \not = 2$, then either $X_1|(k+2)\, x_1$ for $k > 2$ or $X_1| 6\,
x_1$ for $k < 2$, which are impossible. Hence, we should have $k
= 2$.

Let $k = 2$, \eqref{eq:41} becomes
\begin{equation}
\label{eq:42} \begin{array}{rl} &(x_2 - \mu\,
(\frac{x_1^3}{3} - x_1))\, ((x_2 -
\mu\,(\frac{x_1^3}{3} - x_1))\, \dfrac{\partial a_1}{\partial x_1}
-x_1\, \dfrac{\partial a_1}{\partial x_2})\\
 = &(-2\,
\mu\, (x_1^2 - 1)\, (x_2 - \mu (\frac{x_1^3}{3} - x_1)) - 4\,
x_1)\,a_1 - 6\, x_1.
\end{array}
\end{equation}
Let
$$a_1(x_1,x_2) = \sum_{i = 0}^m h_i(x_2)\,  x_1^i,$$
where $h_i(x_2)$ are polynomials and $h_m(x_2) \not=0$. Substitute
$a_1(x_1,x_2)$ into \eqref{eq:42}, and comparing the coefficient
of $x^{m+5}$, we have
$$\frac{1}{9}\,  \mu^2\, m\, h_m(x_2) = \frac{2}{3}\, \mu^2\, h_m(x_2),$$
which implies $m = 6$. Comparing the coefficients of $x_1^{i}\ ( 0\leq
i\leq 10)$, we obtain the equations that are satisfied by $h_i(x_2), i
= 0,\cdots , 6$:
\begin{eqnarray*}
0 &=& x_2\, ( -2\, \mu\,   h_0(x_2) +
x_2\,  h_1(x_2) )\\
0 &=& 6 - 2\, ( -2 + {\mu }^2)\,   h_0(x_2)
+ 2\, x_2^2\,  h_2(x_2) - x_2\,  h_0'(x_2)\\
0 &=& 2\, \mu\, x_2\, h_0(x_2) - ( -4 + {\mu }^2 )\, h_1(x_2) +
2\, \mu\, x_2\, h_2(x_2) +
  3\, x_2^2\,  h_3(x_2)\\
&&{} - \mu\, h_0'(x_2) -
  x_2\,h_1'(x_2)\\
 0 &=& \frac{8\, \mu^2}{3}\,  h_0(x_2) + \frac{4\, \mu\, x_2 }{3}
\,h_1(x_2) + 4\,  h_2(x_2) +
  4\, \mu\, x_2\,   h_3(x_2) + 4\, x_2^2\,  h_4(x_2)\\
&&{} - \mu\,   h_1'(x_2) -
  x_2\,  h_2'(x_2)\\
0 &=& 2\, {\mu }^2\,  h_1(x_2) + \frac{2\, \mu\, x_2 }{3}\,
h_2(x_2) + 4\, h_3(x_2) +
  {\mu }^2\,  h_3(x_2) + 6\, \mu\, x_2\,   h_4(x_2)\\
  &&{} + 5\, x_2^2\,  h_5(x_2) +
  \frac{\mu }{3}\,  h_0'(x_2) - \mu\,   h_2'(x_2) -
  x_2\,  h_3'(x_2)\\
0 &=& \frac{1}{3}(-2\, {\mu }^2 \, h_0(x_2) + 4\, {\mu }^2\,
h_2(x_2) + 12\, h_4(x_2) +
    6\, {\mu }^2\,  h_4(x_2) + 24\, \mu\, x_2\,   h_5(x_2)\\
    &&{} + 18\, x_2^2\,  h_6(x_2) +
    \mu\,   h_1'(x_2) - 3\, \mu\, h_3'(x_2) -
    3\, x_2\,  h_4'(x_2))
\end{eqnarray*}
 \begin{eqnarray*}
0 &=& \frac{1}{9}\,(-5\, {\mu }^2\,  h_1(x_2) + 6\, {\mu }^2\,
h_3(x_2) - 6\, \mu\, x_2\,   h_4(x_2) +
    36\,  h_5(x_2) + 27\, {\mu }^2\,  h_5(x_2)\\
    &&{} + 90\, \mu x_2 \,  h_6(x_2) +
    3\, \mu \,  h_2'(x_2) - 9\, \mu\,   h_4'(x_2) -
    9\, x_2\,  h_5'(x_2))\\
0 &=& -\frac{4\,{\mu }^2}{9}\,  h_2(x_2) - \frac{4\, \mu\, x_2
}{3}\, h_5(x_2) + 4 \, h_6(x_2) +
  4\, {\mu }^2\,  h_6(x_2) + \frac{\mu }{3} \, h_3'(x_2)\\
  &&{} - \mu\,   h_5'(x_2) -
  x_2\,  h_6'(x_2)\\
0 &=& -\frac{\mu}{3}\,( \mu \,  h_3(x_2) + 2\, \mu\,   h_5(x_2) +
6\, x_2\, h_6(x_2) -
         h_4'(x_2) + 3\,  h_6'(x_2) )\\
0 &=& -\frac{\mu}{9}\,( 2\, \mu\,   h_4(x_2) + 12\, \mu\, h_6(x_2)
-
3\,  h_5'(x_2) )\\
0 &=&- \frac{\mu}{9}\,( \mu\,   h_5(x_2) - 3 \, h_6'(x_2))
\end{eqnarray*}
The above equations can be reduced to
$$x_2\, (3\, x_2\, h_5'(x_2) - 2\, \mu\, h_4'(x_2)) = 2\, \mu^3,$$
which is impossible since $h_4(x_2)$
and $h_5(x_2)$ are polynomials. The contradiction concludes that \eqref{eq:37} has
no rational function solution, and hence the order of the
differential Galois group of the van der Pol equation is infinity.
\end{proof}

\section{Proof of Theorem \ref{th:8}}
\label{sec:app} 
Before proving Theorem \ref{th:8}, we introduce some 
notations as following. Let $\delta_1 = \frac{\partial\ }{\partial
x_1},\ \delta_2 = \frac{\partial\ }{\partial x_2}$, $y$ be an
indeterminate over $K$, and denote $\delta_2^iy$ by $y_i$ ($y_0 = y$).  Let
$$X = X_1\,\delta_1 + X_2\,\delta_2,\ \ \ \ \delta_2 X = (\delta_2X_1)\,\delta_1 + (\delta_2X_2)\,\delta_2,$$
$$\hat{X} = X_1\delta_1 + X_2\delta_2 + \sum_{i\geq0}X(y_i)\frac{\partial\ }{\partial y_i},$$
$$B_0 = -X_1\,\delta_2(\frac{X_2}{X_1}),\ \ \ B_i = X_1\,\delta_2(\frac{B_{i-1}}{X_1}) = -X_1\,\delta_2^{i+1}(\frac{X_2}{X_1}),\ \  i = 1,2,\cdots$$
For $F(y)\in K\{\A_0(y)\}$, and $\{X(y)\}$ be the differential ideal that is generated by $X(y)$, we write
\begin{equation}
F(y)\sim R(y)
\end{equation}
if $R(y)\in K\{\A_0(y)\}$ such that $F(y) - R(y)$ is contained in $\{X(y)\}$.

Let $\mathbf{n} = (n_1,n_2,\cdots, n_r)\in {\mathbb{Z}^*}^r$,
define the operators $d_i^j$ and $b_i^j$ for $1<i<j$ by
\begin{equation}
d_i^j(\mathbf{n}) = (n_1,\cdots, n_{j-i}+1,\cdots,n_j-1,\cdots,n_r)
\end{equation}
and
\begin{equation}
b_i^j(\mathbf{n}) = (n_1,\cdots,n_{j-i}-1,\cdots,n_j+1,\cdots,n_r),
\end{equation}
respectively.
Then $b_{i}^j(d_i^j(\mathbf{n})) = \mathbf{n}$. 

Let $\mathbf{n},
\mathbf{m}\in {\mathbb{Z}^*}^r$, the \textit{degree} of
$\mathbf{n}$ is higher than that of $\mathbf{m}$, denoted by
$\mathbf{n}
> \mathbf{m}$, if there exists $1\leq k\leq r$ such that $n_k > m_k$ and
$$n_i = m_i,\ \ i = k+1, \cdots, r.$$
We say $\mathbf{n} \succ \mathbf{m}$ if there exist $1 < i < j$
such that
$$d_i^j(\mathbf{n}) = \mathbf{m}.$$
It is obvious that when $1 > i> j$,
\begin{equation}
\label{eq:a.11} b_i^j(\mathbf{n}) > \mathbf{n} > d_i^j(\mathbf{n})
\end{equation}
and
\begin{equation}
\label{eq:a.12} d_i^j(\mathbf{n}) \succ \mathbf{n} \succ
b_i^j(\mathbf{n}).
\end{equation}

Let $\Lambda$ be the regular prime ideal of QDP 
corresponding to the differential Galois group in Theorem
\ref{th:3}, and assume that $\mathrm{ord}(\Lambda) = r$. Let $A$
in $\Lambda$ with the lowest rank and therefore irreducible. Denoted $A$ as
$$A(x_1,x_2,y,y_1,\cdots, y_r) = \sum_{\mathbf{m}}A_{\mathbf{m}}(x_1,x_2,y)y_1^{m_1}\,\cdots\, y_r^{m_r}$$
and let
$$\mathcal{I}_A = \{\mathbf{m}\in {\mathbb{Z}^*}^r|\ \ A_{\mathbf{m}}\not=0\}.$$
By $\mathbf{n}$ we will always denote the element in $\mathcal{I}_A$ with the highest
degree. We assume further that that $A_{\mathbf{n}} = 1$ and
the coefficients $A_{\mathbf{m}}$ are rational functions
in $K(\A_0(y))$. For any $\mathbf{m}\in \mathcal{I}_A$, let
\begin{equation}
\label{eq:a.13} \mathcal{P}(\mathbf{m}) = \{\mathbf{p}\in \mathcal{I}_A\ |\
d_i^j(\mathbf{p})\succ \mathbf{m}\ \mathrm{for\ some}\ 1< i<j\}
\end{equation}
 and
$\#(\mathbf{m}) = |\mathcal{P}(\mathbf{m})|$. A subset $\mathcal{J}_A\subseteq \mathcal{I}_A$ is
\textit{closed} if for every $\mathbf{m}\in \mathcal{J}_A$, $\mathbf{p}\succ
\mathbf{m}$ implies $\mathbf{p}\in \mathcal{J}_A$. It is easy to see that $\mathcal{I}_A$
and $\{\mathbf{n}\}$ are closed.

The proof will be done by showing that all possible
structures of $\mathcal{I}_A$ are, besides the cases with $r = 0$ and $r =
\infty$,
\begin{enumerate}
\item $r = 1$, and $\mathcal{I}_A = \{n,0\}$; or
\item $r = 2$, and
$\mathcal{I}_A = \{(0,1), (1,0)\}$, with
$$(0,1) \succ (1,0);$$
or
\item $r = 3$, and $\mathcal{I}_A = \{(1,0,1), (0,2,0), (2,0,0)\}$, with relations
\begin{center}
\unitlength=0.5cm
\begin{picture}(8,4.5)
\put(2,0){(2,0,0)} \put(2,2){(1,1,0)} \put(2,4){(1,0,1)}
\put(5,2){(0,2,0)} \put(2.8,1){$\curlyvee$}
\put(2.8,3){$\curlyvee$} \put(4.2,2){$\prec$}
\put(1.0,4.2){\line(1,0){0.8}} \put(1.0,4.2){\line(0,-1){1.2}}
\put(1.0,0.2){\line(0,1){1.2}}\put(1.0,0.2){\line(1,0){0.8}}
\put(0.8,2.0){$\curlyvee$}
\end{picture}
\end{center}
Here $(1,1,0)$ is an auxiliary index with $A_{(1,1,0)} = 0$.
\end{enumerate}

The proof will be complete following the flow chart in Figure \ref{fig:1}.

\begin{figure}[hbt]
\centering
\unitlength=0.5cm
\begin{picture}(22,9)
\put(18,8){\framebox{Lemma \ref{lem:a.1}}}
\put(12,8){\framebox{Lemma \ref{lem:a.2}}}
\put(6,8){\framebox{Lemma \ref{lem:a.6}}}
\put(0,8){\framebox{Lemma \ref{lem:a.9}}}
\put(0,6){\framebox{Lemma \ref{lem:a.8}}}
\put(0,4){\framebox{Lemma \ref{lem:a.12}}}
\put(6,4){\framebox{Lemma \ref{lem:a.10}}} \put
(12,6){\framebox{Lemma \ref{lem:a.3}}} \put(18,6){\framebox{Lemma
\ref{lem:a.4}}} \put(18,4){\framebox{Lemma \ref{lem:a.15}}}
\put(0,2){\framebox[11cm]{Lemma \ref{lem:a.7}}}
\put(0,0){\framebox{Lemma \ref{le:2}}} \put(6,6){\framebox{Lemma
\ref{lem:a.11}}} \put(6,0){\framebox[8cm]{Theorem \ref{th:8}}}
\put(18,8.2){\vector(-1,0){1.5}} \put(12,8.2){\vector(-1,0){1.5}}
\put(6,8.2){\vector(-1,0){1.5}} \put(2,7.8){\vector(0,-1){1.0}}
\put(4.5,6.2){\vector(1,0){1.3}} \put(8,5.8){\vector(0,-1){1.0}}
\put(4.5,4.2){\vector(1,0){1.5}} \put(4.5,6){\vector(1,-1){1.5}}
\put(16.5,6.2){\vector(1,0){1.5}} \put(12,6.2){\vector(-1,0){1.1}}
\put(18,6){\vector(-1,-2){1.6}} \put(20,3.8){\vector(0,-1){1.0}}
\put(14,5.8){\vector(0,-1){3.0}}  \put(14,1.7){\vector(0,-1){1.0}}
\put(2,3.8){\vector(0,-1){1.0}} \put(8,3.8){\vector(0,-1){1.0}}
\put(10.5,5.8){\vector(1,-2){1.5}} \put(4.5,0.2){\vector(1,0){1.5}}
\put(0,6.2){\line(-1,0){0.5}} \put(-0.5,6.2){\line(0,-1){4.0}}
\put(-0.5,2.2){\vector(1,0){0.5}} \put(20,6.7){\line(0,1){0.6}}
\put(20,7.3){\line(-1,0){12.0}}\put(8,7.3){\vector(0,-1){0.6}}
\end{picture}
\caption{Flow chart of the proof of Theorem \ref{th:8}}
\label{fig:1}
\end{figure}

\begin{lem}
\label{lem:a.3}If $u\not=0$ satisfies
\begin{equation}
\label{eq:a.8} X\,u = B_0\,u
\end{equation}
then there exists a first integral $\omega$ of $X$ such that
$$\delta_2\omega = u.$$
\end{lem}
\begin{proof}From \eqref{eq:a.8}, we have
\begin{eqnarray*}
X_1\,\delta_1u + X_2\,\delta_2 u &=& B_0\,u =
-X_1\,\delta_2(\frac{X_2}{X_1})\,u\\
\delta_1\,u + \frac{X_2}{X_1}\,\delta_2u &=&
-\delta_2(\frac{X_2}{X_1}\,u)\\
\delta_1u&=&-\delta_2(\frac{X_2}{X_1}\,u -
\frac{X_2}{X_1}\,\delta_2u)\\
&=&\delta_2(-\frac{X_2}{X_1}\,u)
\end{eqnarray*}
Let $v =-\frac{X_2}{X_1}\,u$, then the 1-form $v dx_1  + u dx_2$
is closed, and
$$\omega(x_1,x_2) = \int_{(0,0)}^{(x_1,x_2)} v dx_1  + u dx_2$$
is a first integral of $X$, with $\delta_2\omega = u$.
\end{proof}

\begin{lem}
\label{lem:a.4}If there exists $u$ satisfying
\begin{equation}
\label{eq:a.9} X u = B_0\, u  + B_1
\end{equation}
 then $X$ has a first
integral $\omega$ such that
$$\frac{\delta_2^2\omega}{\delta_2 \omega} = u.$$
\end{lem}
\begin{proof}From \eqref{eq:a.9}, we have
\begin{eqnarray*}
X_1\,\delta_1 u + X_2\,\delta_2 u &=&
-X_1\,\delta_2(\frac{X_2}{X_1})\,u +
X_1\,\delta_2(\frac{B_0}{X_1})\\
\delta_1 u&=& -\frac{X_2}{X_1}\,\delta_2u -
\delta_2(\frac{X_2}{X_1})\,u + \delta_2(\frac{B_0}{X_1})\\
&=&\delta_2(-\frac{X_2}{X_1}\,u + \frac{B_0}{X_1}).
\end{eqnarray*}
Thus, let
$$v = -\frac{X_2}{X_1}\,u + \frac{B_0}{X_1},$$
the 1-form $v\, d x_1 + u\, d x_2$ is a closed. Let
$$\eta(x_1,x_2) = \exp\left[\int_{(0,0)}^{(x_1,x_2)} v\, d x_1 + u\, d x_2\right],$$
then
$$X(\eta) = \eta\,(X_1 v + X_2 u) = \eta\, (X_1 (-\frac{X_2}{X_1}\,u + \frac{B_0}{X_1}) + X_2 u) = B_0\,\eta.$$
From Lemma \ref{lem:a.3}, there exists a first integral $\omega$ of
$X$ such that
$$\delta_2\omega  = \eta,$$
and therewith
$$\frac{\delta_2^2\omega}{\delta_2\omega} = u.$$
The Lemma is concluded.
\end{proof}

\begin{lem}
\label{lem:a.15}If there exists $u$ satisfying
\begin{equation} \label{eq:a.10} X u = 2\,B_0\, u  + B_2,
\end{equation}
then $X$ has a first integral $\omega$ of $X$ such that
$$\frac{2\,\delta_2\omega\cdot \delta_2^3\omega - 3(\delta_2^2\omega)^2}{(\delta_2 \omega)^2} = u.$$
\end{lem}
\begin{proof}From \eqref{eq:a.10}, we have
\begin{eqnarray*}
X_1\,\delta_1 u + X_2\,\delta_2 u &=&
-2\,X_1\,\delta_2(\frac{X_2}{X_1})u -
X_1\,\delta_2^3(\frac{X_2}{X_1}),\\
\delta_1u + \frac{X_2}{X_1}\delta_2u&=&-2\,\delta_2(\frac{X_2}{X_1})u - \delta_2^3(\frac{X_2}{X_1})\\
\delta_1 u&=&-2\,\delta_2(\frac{X_2}{X_1})\,u -
\frac{X_2}{X_1}\,\delta_2u - \delta_2^3(\frac{X_2}{X_1}),\\
&=&-\delta_2(\frac{X_2}{X_1}\,u) - \frac{X_2}{X_1}\,\delta_2u -
\delta_2^3(\frac{X_2}{X_1}).
\end{eqnarray*}
Consider the partial differential equations:
\begin{equation}
\label{eq:a.4}\left\{
\begin{array}{rcl}
\delta_2w &=& u + \frac{1}{2}\,w^2\\
\delta_1w &=& -\delta_2^2(\frac{X_2}{X_1})-\frac{X_2}{X_1}\,u -
\delta_2(\frac{X_2}{X_1})\,w - \frac{1}{2}(\frac{X_2}{X_1})\,w^2.
\end{array}\right.
\end{equation}
We have
\begin{eqnarray*}
\delta_1\delta_2w&=&\delta_1(u + \frac{1}{2}w^2)\\
&=&\delta_1u + w\,\delta_1w\\
&=&-\delta_2(\frac{X_2}{X_1}u) - \frac{X_2}{X_1}\,\delta_2u -
\delta_2^3(\frac{X_2}{X_1})\\
&&{} +
w\,\left(-\delta_2^2(\frac{X_2}{X_1})-\frac{X_2}{X_1}\,u -
\delta_2(\frac{X_2}{X_1})\,w -
\frac{1}{2}(\frac{X_2}{X_1})\,w^2\right)\\
&=&-\delta_2(\frac{X_2}{X_1}u) - \frac{X_2}{X_1}\,\delta_2u -
\delta_2^3(\frac{X_2}{X_1}) -
\left(\delta_2^2(\frac{X_2}{X_1})+\frac{X_2}{X_1}\,u\right)\,w\\
&&{} -
\delta_2(\frac{X_2}{X_1})\,w^2 -
\frac{1}{2}(\frac{X_2}{X_1})\,w^3
\end{eqnarray*}
and
\begin{eqnarray*}
\delta_2\delta_1w&=&\delta_2\left(-\delta_2^2(\frac{X_2}{X_1})-\frac{X_2}{X_1}\,u
- \delta_2(\frac{X_2}{X_1})\,w -
\frac{1}{2}(\frac{X_2}{X_1})\,w^2\right)\\
&=&-\delta_2^3(\frac{X_2}{X_1}) - \delta_2(\frac{X_2}{X_1}u) -
\delta_2^2(\frac{X_2}{X_1})w -
\delta_2(\frac{X_2}{X_1})\,\delta_2w\\
&&{} -
\frac{1}{2}\delta_2(\frac{X_2}{X_1})\,w^2 -
\frac{X_2}{X_1}\,w\,\delta_2w\\
&=&-\delta_2^3(\frac{X_2}{X_1}) - \delta_2(\frac{X_2}{X_1}u) -
\delta_2^2(\frac{X_2}{X_1})w - \delta_2(\frac{X_2}{X_1})\,(u +
\frac{1}{2}\,w^2)\\
&&{} - \frac{1}{2}\delta_2(\frac{X_2}{X_1})\,w^2 -
\frac{X_2}{X_1}\,w\,(u + \frac{1}{2}\,w^2)\\
&=&-\delta_2^3(\frac{X_2}{X_1}) -
\delta_2(\frac{X_2}{X_1}u)-\delta_2(\frac{X_2}{X_1})\,u -
\left(\delta_2^2(\frac{X_2}{X_1}) + \frac{X_2}{X_1}u\right)\,w\\
&&{} -
\delta_2(\frac{X_2}{X_1})\,w^2 - \frac{1}{2}(\frac{X_2}{X_1})\,w^3
\end{eqnarray*}
Therefore, $\delta_1\delta_2w = \delta_2\delta_1w$, and the
equations \eqref{eq:a.4} have a solution $w$ that is analytic at
$(0,0)$. Let
$$v = -\delta_2(\frac{X_2}{X_1}) - (\frac{X_2}{X_1})\,w, $$
then
\begin{eqnarray*}
\delta_2v&=&\delta_2(-\delta_2(\frac{X_2}{X_1}) -
(\frac{X_2}{X_1})\,w)\\
&=&-\delta_2^2(\frac{X_2}{X_1}) -
\delta_2(\frac{X_2}{X_1})\,w - (\frac{X_2}{X_1})\,\delta_2w\\
&=&-\delta_2^2(\frac{X_2}{X_1}) -
\delta_2(\frac{X_2}{X_1})\,w - (\frac{X_2}{X_1})\,(u + \frac{1}{2}\,w^2)\\
&=&-\delta_2^2(\frac{X_2}{X_1}) - (\frac{X_2}{X_1})\,u -
\delta_2(\frac{X_2}{X_1})w - \frac{1}{2}(\frac{X_2}{X_1})\,w^2\\
&=&\delta_1w.
\end{eqnarray*}
Therefore, the 1-form $v dx_1 + w dx_2$ is a closed. Let
$$\omega_2 =
\exp\left[\int_{(0,0)}^{(x_1,x_2)} v dx_1 + w dx_2\right],\ \ \
\omega_1 = -\frac{X_2}{X_1}\omega_2,$$ than
\begin{eqnarray*}
\delta_1\omega_2&=&\omega_2\,v\\
\delta_2\omega_1&=&-\delta_2(\frac{X_2}{X_1})\,\omega_2 -
(\frac{X_2}{X_1})\,\delta_2\omega_2\\
&=&\omega_2(-\delta_2(\frac{X_2}{X_1}) -
\frac{X_2}{X_1}\,\frac{\delta_2\omega_2}{\omega_2})\\
&=&\omega_2(-\delta_2(\frac{X_2}{X_1}) - \frac{X_2}{X_1}\,w)\\
&=&\omega_2\,v = \delta_1\,\omega_2.
\end{eqnarray*}
Hence, $\omega = \int_{(0,0)}^{(x_1,x_2)} \omega_1 dx_1 + \omega_2
dx_2$ is well defined, and is a first integral of $X$ at $(0,0)$.
It is easy to verify that
$$\dfrac{2\delta_2\omega\cdot\delta_2^3\omega - 3(\delta_2^2\omega)^2}{(\delta_2\omega)^2} = u.$$
The proof is completed.
\end{proof}

\begin{lem}
\label{lem:a.1}We have
\begin{enumerate} 
\item[(1).] $\delta_2X =
(\frac{\delta_2 X_1}{X_1})\,X - B_0\,\delta_2$; 
\item[(2).] $X(y_j) =
\delta_2(X(y_{j-1})) - (\frac{\delta_2 X_1}{X_1})\,X(y_{j-1}) +
B_0\,y_j$.
\end{enumerate}
\end{lem}
\begin{proof}
The proof is straightforward from (1)
\begin{eqnarray*} \delta_2X &=& (\delta_2X_1)\,\delta_1 +
(\delta_2X_2)\,\delta_2\\
&=&\frac{\delta_2 X_1}{X_1}\, (X_1\,\delta_1 + X_2\,\delta_2) -
\frac{X_2}{X_1} (\delta_2 X_1)\,\delta_2 + (\delta_2
X_2)\,\delta_2\\
&=&\frac{\delta_2 X_1}{X_1}\,X + X_1\,\frac{X_1\,\delta_2X_2 -
X_2\,\delta_2 X_1}{X_1^2}\delta_2\\
&=&\frac{\delta_2X_1}{X_1}\,X - B_0\,\delta_2,
\end{eqnarray*}
and (2)
\begin{eqnarray*}
X(y_j) &=& \delta_2(X(y_{j-1})) - (\delta_2 X)y_{j-1}\\
&=& \delta_2(X(y_{j-1})) - (\frac{\delta_2X_1}{X_1}\,X -
B_0\,\delta_2)(y_{j-1})\\
&=& \delta_2(X(y_{j-1})) - \frac{\delta_2X_1}{X_1}\,X(y_{j-1}) +
B_0\,\delta_2(y_{j-1})\\
&=&\delta_2(X(y_{j-1})) - (\frac{\delta_2 X_1}{X_1})\,X(y_{j-1}) +
B_0\,y_j.
\end{eqnarray*}
\end{proof}

\begin{lem}
\label{lem:a.2} We have
\begin{equation} \label{eq:a.14} X(y_j) \sim
\sum_{i = 0}^{j-1}a_{j,i}\,B_i\,y_{j-i}
\end{equation}
 where $a_{j,i}$ are
constants, with $a_{j,0} = j$.
\end{lem}
\begin{proof} From Lemma \ref{lem:a.1}, when $j = 1$, we have
$$
X(y_1) = \delta_2(X(y_0)) - (\frac{\delta_2 X_1}{X_1})\,X(y_0) +
B_0\,y_1 \sim B_0\,y_1,
$$
which is the desired \eqref{eq:a.14} with $a_{1,0} = 1$. 

Assume that \eqref{eq:a.14} is valid for $j = k$, and $a_{k,0}=k$, then
by Lemma \ref{lem:a.1},
\begin{eqnarray*}
X(y_{k+1}) &=& \delta_2(X(y_k)) - (\frac{\delta_2
X_1}{X_1})\,X(y_k) + B_0\,y_{k+1}\\
&\sim&\delta_2(\sum_{i=0}^{k-1}a_{k,i}B_i y_{k-i}) -
(\frac{\delta_2 X_1}{X_1})\,(\sum_{i=0}^{k-1}a_{k,i}B_i y_{k-i})
+ B_0\,y_{k+1}\\
&=&\sum_{i=0}^{k-1}a_{k,i}\,((\delta_2B_i)\,y_{k-i} +
B_i\,\delta_2y_{k-i}) -
\sum_{i=0}^{k-1}a_{k,i}\frac{\delta_2\,X_1}{X_1}\,B_i\,y_{k-i} +
B_0y_{k+1}\\
&=&\sum_{i=0}^{k-1}a_{k,i}\left((\delta_2B_i - \frac{\delta_2
X_1}{X_1}\,B_i)\,y_{k-i} + B_i\,y_{k-i+1}\right) + B_0\,y_{k+1}\\
&=&(a_{k,0} + 1)\,B_0\,y_{k+1} + \sum_{i = 0}^{k-2}(a_{k,i}
X_1\,\delta_2(\frac{B_i}{X_1}) + a_{k,i+1}\,B_{i+1})y_{k-i}\\
&&{} + a_{k,k-1}\,X_1\delta_2(\frac{B_{k-1}}{X_1})\,y_1\\
&=&(a_{k,0} + 1)\,B_0\,y_{k+1} + \sum_{i = 0}^{k-2}(a_{k,i}
 + a_{k,i+1})\,B_{i+1}y_{k-i} + a_{k,k-1}\,B_{k}\,y_1\\
 &=&\sum_{i = 0}^ka_{k+1,i}\,B_i\,y_{k + 1 - i},
\end{eqnarray*}
where
$$\left\{\begin{array}{ll}
a_{k+1,0} = a_{k,0} + 1 = k+1,& \\
a_{k+1,i} = a_{k,i-1} + a_{k,i},&\ \ (1\leq i\leq k-1),\\
a_{k+1,k} = a_{k,k-1}.&
\end{array}\right.$$
Thus, the Lemma has been proved.
\end{proof}

\begin{lem}
\label{lem:a.12} Let $\mathbf{m}\in {\mathbb{Z}^*}^r$, and define
\begin{equation}
\label{eq:a.5} C({\mathbf{m}}) = m_1 + 2\,m_2 + \cdots + r\,m_r.
\end{equation}
 If
$\mathbf{m}\succ \mathbf{p}$, then $C({\mathbf{m}}) >
C({\mathbf{p}})$. In particularly, if $d_i^j(\mathbf{m}) =
\mathbf{p}$, then $C(\mathbf{m}) - C({\mathbf{p}}) = i$.
\end{lem}
\begin{proof} Let $d_i^j(\mathbf{m}) = \mathbf{p}$, then
$$C({\mathbf{m}}) - C({\mathbf{p}}) = (j-i)\,m_{j-i} + j\,m_j - ((j-i)\,(m_{j-i}+1) + j\,(m_j-1)) = i.$$
\end{proof}

\begin{lem}
\label{lem:a.6} Let $P = A_{\mathbf{m}}y_1^{m_1}\cdots y_r^{m_r}$,
then
$$
\hat{X}(P)\sim (X(A_{\mathbf{m}}) + C({\mathbf{m}})\, B_0\,
A_{\mathbf{m}})\,\mathbf{y}^{\mathbf{m}} + \sum_{i =
1}^{r-1}\sum_{j=i+1}^rn_j\,a_{j,i}\,B_i\,A_\mathbf{m}\,\mathbf{y}^{d_i^j(\mathbf{m})},
$$
where $\mathbf{y}^{\mathbf{m}} = y^{m_1}\cdots y_r^{m_r}$.
\end{lem}
\begin{proof} By Lemma \ref{lem:a.2}, we have
\begin{eqnarray*}
\hat{X}(P) &=& X(A_\mathbf{m})\mathbf{y}^{\mathbf{m}} +
\frac{\partial A_{\mathbf{m}}}{\partial
y}\,\mathbf{y}^{\mathbf{m}}\,X(y) + A_\mathbf{m}\,\sum_{j=1}^r m_j
y_1^{m_1}\cdots
y_j^{m_j-1}\cdots y_r^{m_r}\,X(y_i)\\
&\sim&X(A_\mathbf{m})\mathbf{y}^{\mathbf{m}} +
A_\mathbf{m}\,\sum_{j=1}^r m_j y_1^{m_1}\cdots y_j^{m_j-1}\cdots
y_r^{m_r}\,(\sum_{i=0}^{j-1}a_{j,i}\,B_i\,y_{j-i})\\
&=&X(A_\mathbf{m})\mathbf{y}^{\mathbf{m}} +
A_\mathbf{m}\,B_0(\sum_{j =
1}^rm_j\,a_{j,0})\mathbf{y}^{\mathbf{m}}\\
&&{} + A_\mathbf{m}\sum_{j=1}^r\sum_{i = 1}^{j-1}m_ja_{j,i}B_i\,
y_1^{m_1}\cdots
y_{j-i}^{m_{j-i}+1}\,\cdots\,y_j^{m_j-1}\cdots y_r^{m_r}\\
&=&(X(A_{\mathbf{m}}) + C({\mathbf{m}})\, B_0\,
A_{\mathbf{m}})\,\mathbf{y}^{\mathbf{m}} + \sum_{i =
1}^{r-1}\sum_{j=i+1}^rm_j\,a_{j,i}\,B_i\,A_\mathbf{m}\,\mathbf{y}^{d_i^j(\mathbf{m})}
\end{eqnarray*}
and the Lemma is concluded.
\end{proof}

\begin{lem}
\label{lem:a.9} Let $A\in \Lambda$ with the lowest rank and $r =
\mathrm{ord}(\Lambda) > 1$. Let $\mathbf{n}\in \mathcal{I}_A$ with the
highest degree and $A_\mathbf{n} = 1$, then for any $\mathbf{m} <
\mathbf{n}$,
\begin{equation}
\label{eq:a.2} X(A_\mathbf{m}) = (C({\mathbf{n}}) -
C({\mathbf{m}}))\,B_0\,A_{\mathbf{m}}
-\sum_{i=1}^{r-1}\sum_{j=i+1}^{r}(m_j+1)\,a_{j,i}\,B_i\,A_{{b_i^j}(\mathbf{m})}
\end{equation}
where $A_{\mathbf{m}} = 0$ if $\mathbf{m}\not\in \mathcal{I}_A$.
\end{lem}
\begin{proof}
Let
$$A = \sum_{\mathbf{m}\in \mathcal{I}_A}A_{\mathbf{m}}\mathbf{y}^{\mathbf{m}}.$$
By Lemma \ref{lem:a.6}, we have
\begin{eqnarray*}
X(A)&=&\sum_{\mathbf{m}\in \mathcal{I}_A}X(A_{\mathbf{m}}\mathbf{y}^{\mathbf{m}})\\
&\sim&\sum_{\mathbf{m}\in \mathcal{I}_A}\left((X(A_{\mathbf{m}}) +
C({\mathbf{m}})\,B_0\,A_{\mathbf{m}})\,\mathbf{y}^{\mathbf{m}}  +
\sum_{i=1}^{r-1}\sum_{j =
i+1}^{r}m_j\,a_{j,i}B_i\,A_{\mathbf{m}}\mathbf{y}^{d_{i}^j(\mathbf{m})}\right)\\
&=&\sum_{\mathbf{m}\in \mathcal{I}_A}\left(X(A_{\mathbf{m}}) +
C({\mathbf{m}})\,B_0\,A_{\mathbf{m}} + \sum_{i=1}^{r-1}\sum_{j =
i+1}^{r}(m_j+1)\,a_{j,i}B_i\,A_{{b_i^j}(\mathbf{m})}\right)\,\mathbf{y}^{\mathbf{m}}
\end{eqnarray*}
Note that for any $j>i$, $b_i^j(\mathbf{n}) > \mathbf{n}$,  and
thus $b_i^j(\mathbf{n}) \not\in \mathcal{I}_A$, i.e., $A_{b_i^j(\mathbf{n})}
= 0$. Taking account $A_{\mathbf{n}} = 1$, we have
\begin{eqnarray}
\label{eq:l58}
X(A) - C({\mathbf{n}})\,B_0\,A  &=& \sum_{\mathbf{m} < \mathbf{n}}\Big(X(A_\mathbf{m})
+ (C({\mathbf{m}}) - C({\mathbf{n}}))\,B_0\,A_{\mathbf{m}}\nonumber\\
&&{} +
\sum_{i=1}^{r-1}\sum_{j=i+1}^{r-1}(m_j+1)\,a_{j,i}\,B_i\,A_{{b_i^j}(\mathbf{m})}\Big)\mathbf{y}^{\mathbf{m}}
\end{eqnarray}
Thus, $X(A) - C({\mathbf{n}})\,B_0\,A$ contains in
$\Lambda$ and has lower rank than $A$. But $A$ is the element in
$\Lambda$ with the lowest rank, thus  
$X(A)- C({\mathbf{n}})\,B_0\,A\equiv0$. Therefore the coefficients in \eqref{eq:l58} are zero, i.e.
$$X(A_{\mathbf{m}}) +
(C({\mathbf{m}}) - C(\mathbf{n}))\,B_0\,A_{\mathbf{m}} + \sum_{i=1}^{r-1}\sum_{j =
i+1}^{r-1}(m_j+1)\,a_{j,i}B_i\,A_{{b_i^j}(\mathbf{m})} = 0,\ \ (\forall \mathbf{m}<\mathrm{n})$$
from which \eqref{eq:a.2} is concluded.
\end{proof}

From Lemma \ref{lem:a.9}, we have
\begin{lem}
\label{lem:a.8} If $\mathcal{P}(\mathbf{m}) = \{\mathbf{p}_1,\cdots,
\mathbf{p}_k\}$, and $d_{i_l}^{j_l}(\mathbf{p}_l) = \mathbf{m},\
(l = 1,2,\cdots,k)$, then the coefficients $A_{\mathbf{p}_l},
A_{\mathbf{m}}$ satisfy
\begin{equation}
\label{eq:a.3} X(A_{\mathbf{m}}) = (C({\mathbf{n}}) -
C({\mathbf{m}}))\,B_0\,A_{\mathbf{m}} -\sum_{l = 1}^k
(m_{j_l}+1)\,a_{j_l,i_l}\,B_{i_l}\,A_{\mathbf{p}_l}.
\end{equation}
\end{lem}

\begin{lem}
\label{lem:a.10}Let $A\in \Lambda$ with the lowest rank and $r =
\mathrm{ord}(\Lambda) > 1$. Let $\mathbf{n}\in \mathcal{I}_A$ with the highest
degree. Then for any $\mathbf{m}\in \mathcal{I}_A$, $\#(\mathbf{m})=0$ if
and only if $C(\mathbf{m}) = C(\mathbf{n})$. Furthermore, if
$\#(\mathbf{m})=0$, then $A_{\mathbf{m}}$ is a constant.
\end{lem}
\begin{proof} At first we will prove that if $\#(\mathbf{m}) = 0$, then $C(\mathbf{m}) = C(\mathbf{n})$.

When $\#(\mathbf{m}) = 0$,  by Lemma \ref{lem:a.8}, we have
$$X(A_{\mathbf{m}}) = (C({\mathbf{n}}) -
C({\mathbf{m}}))\,B_0\,A_{\mathbf{m}}.$$ 
If
$C({\mathbf{m}})\not=C({\mathbf{n}})$, let $n = C(\mathbf{n}) -
C(\mathbf{m})$, then
$$X(A_{\mathbf{m}}^{1/n})
 = B_0\,A_{\mathbf{m}}^{1/n}.$$
 From Lemma \ref{lem:a.3}, there exists a first integral
 $\omega$ of $X$ such that
 $$\delta_2\omega = A_{\mathbf{m}}^{1/n},$$
 and hence
 $$y_1^{n} - A_{\mathbf{m}}\in \Lambda,$$
 which yields $r = 1$ and contradict. Thus, we have concluded that $C(\mathbf{m}) =
 C(\mathbf{n})$.

Now, we will prove that if $C(\mathbf{m}) = C(\mathbf{n})$, then
$\#(\mathbf{m}) = 0$. If on the contrary,
 $\#(\mathbf{m})>0$, then there exists $\mathbf{p}\in
 \mathcal{P}(\mathbf{m})$, and from \eqref{eq:a.11}, $C(\mathbf{p}) > C(\mathbf{m}) =
 C(\mathbf{n})$. On the other hand, there exist $\mathbf{p}_1, \cdots, \mathbf{p}_k$, such that
$$\mathbf{p}_1 \succ \cdots \succ \mathbf{p}_k = \mathbf{p}$$
and $\#(\mathbf{p}_1) = 0$. Therefore $C(\mathbf{p}) = C(\mathbf{n})$. It is easy to have $C(\mathbf{p})
 = C(\mathbf{p}_k) <\cdots < C(\mathbf{p}_{1}) = C(\mathbf{n})$,
which is contradict, and the statement is concluded.

Now, we have  proved that $\#(\mathbf{m})$ if and only if $C(\mathbf{m}) = C(\mathbf{n})$.

If $\#(\mathbf{m}) = 0$, then $C(\mathbf{n}) = C(\mathbf{m})$, and
by \eqref{eq:a.3}, $X(A_{\mathbf{m}}) = 0$. But
$\mathrm{ord}(\Lambda)
> 0$, thus $A_{\mathbf{m}}$ is a constant.
\end{proof}

\begin{lem}
\label{lem:a.11}Let $A\in \Lambda$ with the lowest rank and $r =
\mathrm{ord}(\Lambda) \geq 3$. Assume that $\mathcal{J}_A\subseteq \mathcal{I}_A$ is
closed, and $\mathbf{m} = (m_1,m_2,\cdots,m_r)\in \mathcal{J}_A$ with the
highest degree, then $m_2 = 0$.
\end{lem}
\begin{proof}
If on the contrary, $m_2>0$, let $\mathbf{p} = d_1^2(\mathbf{m}) =
(m_1+1,m_2-1,m_3,\cdots,m_r)$. It is easy to verify 
$\mathcal{P}(\mathbf{p}) = \{\mathbf{m}\}$ by (1) $d_1^2(\mathbf{m}) = \mathcal{P}$,
and (2) if any other $\mathbf{m'}$ such that $d_i^j(\mathbf{m'})
=\mathcal{P}$, then $\mathbf{m'} > \mathbf{m}$. Hence, we have
\begin{equation}
\label{eq:a.15} \left\{
\begin{array}{rcl}
X(A_{\mathbf{m}}) &=& (C(\mathbf{n}) - C(\mathbf{m}))\,B_0\,A_{\mathbf{m}}\\
X(A_{\mathbf{p}}) &=& B_0\,A_{\mathbf{p}} - m_2\,a_{2,1}\,
B_1\,A_{\mathbf{m}}.
\end{array}\right.
\end{equation}
From \eqref{eq:a.15} and Lemmas \ref{lem:a.3} and
\ref{lem:a.4}, either $r = 0$ (if $C_{\mathbf{n}} =
C_{\mathbf{m}}$ and $A_{\mathbf{m}}$ is not a constant), or $r = 1$
(if $C_{\mathbf{n}} \not= C_{\mathbf{m}}$), or $r = 2$ (if
$C(\mathbf{n}) - C(\mathbf{m}) = 0$ and $A_{\mathbf{m}}$ is
constant), contradict with the assumption $r\geq 3$.
\end{proof}

\begin{lem}
\label{lem:a.7} Assume that $K$ contains no first integral of $X$,
and let $r$ to be the order of the differential Galois group of
$X$.
\begin{enumerate}
\item[(1).] If $r = 1$, there exists a first integral $\omega$ of
$X$, and $n\in \mathbb{N}$, such that
$$(\delta_2\omega)^n\in K.$$
 \item[(2).] If $r = 2$, there exists a first integral $\omega$ of $X$, such that
 $$\delta_2^2\omega/\delta_2\omega\in K.$$
\item[(3).] If $r = 3$, there exists a first integral $\omega$ of
$X$, such that
$$\dfrac{2\,\delta_2\omega\cdot\delta_2^3\omega - 3(\delta_2^2\omega)^2}{(\delta_2\omega)^2}\in K.$$
\item[(4).] If $r\not=\infty$, then $r\leq 3$.
\end{enumerate}
\end{lem}
\begin{proof}
Let $\Lambda$ to be the regular prime ideal of QDP 
corresponding to the differential Galois of $X$, and $A\in \Lambda$
with the lowest rank, $\mathbf{n}\in \mathcal{I}_A$ with the highest degree.

(1). If $r = 1$, we can write $A$ as
$$A = y_1^n + A_1\,y_1^{n-1} + \cdots + A_n$$
with $A_i \in K[\A_0(y)]$ and $A_n\not=0$ since $A$ is irreducible. From Lemma \ref{lem:a.8}, we have
$$X(A_n) = n\,B_0\,A_n,$$
i.e.,
$$X(A_n^{1/n}) = B_0\,A_n^{1/n}.$$
By Lemma \ref{lem:a.3} and $A_n\not=0$, there exists a first integral $\omega$ of $X$ such
that
$$\delta_2\omega = A_n^{1/n},$$
i.e.,
$$(\delta_2\omega)^n = A_n\in K[\A_0(y)].$$
But $\delta_2\omega$ is independent to $y$, hence $A_n$ does not involves $y$, and therefore $(\delta_2\omega)^n \in K$.

(2). If $r = 2$, let $\mathbf{n} = (n_1,n_2)$  and $\mathbf{m} =
d_1^2(\mathbf{n}) = (n_1+1,n_2-1)$, then $\mathcal{P}(\mathbf{m}) =
\{\mathbf{n}\}$. Thus, by Lemma \ref{lem:a.8} and Lemma
\ref{lem:a.12}, we have
$$X(A_{\mathbf{m}}) = B_0\,A_{\mathbf{m}} - n_2\,a_{2,1}\,B_1,$$
i.e.,
$$X(-\frac{A_{\mathbf{m}}}{n_2\,a_{2,1}}) = B_0\,(-\frac{A_{\mathbf{m}}}{-n_2\,a_{2,1}})+B_1.$$
From Lemma \ref{lem:a.4}, there exists a first integral $\omega$
such that
$$\frac{\delta_2^2\omega}{\delta_2\omega} = -\frac{A_{\mathbf{m}}}{n_2\,a_{2,1}}\in K[\A_0(y)].$$
Similar to the above argument, we have $\delta_2^2\omega/\delta_2\omega\in K$ and (2) is proved.

(3). If $r = 3$, we can write $\mathbf{n} = (n_1,n_2,n_3)$. From Lemma
\ref{lem:a.11}, we have $\mathbf{n} = (n_1,0,n_3)$. Let
$$\mathbf{p} = d_1^3(\mathbf{n}) = (n_1,1,n_3-1),$$
$$ \mathbf{q} = d_2^3(\mathbf{n}) = (n_1+1,0,n_3-1),$$
$$ \mathbf{m} = {b_1^2}(\mathbf{p}) = (n_1-1,2,n_3-1).$$ It is easy to have
$C(\mathbf{m}) = C(\mathbf{n})$. Therefore, by Lemma
\ref{lem:a.10}, $A_{\mathbf{m}}$ is a constant. Furthermore, we
have $\mathcal{P}(\mathbf{p}) = \{\mathbf{n},\mathbf{m}\}$ and
$\mathcal{P}(\mathbf{q}) = \{\mathbf{n},\mathbf{p}\}$.  Lemma
\ref{lem:a.8}, yields
\begin{equation}
\label{eq:a.6} X(A_{\mathbf{p}}) = B_0\,A_{\mathbf{p}} -
(n_3\,a_{3,1}\,A_{\mathbf{n}} +
2\,a_{2,1}\,A_{\mathbf{m}})\,B_1,
\end{equation}
and
\begin{equation}
\label{eq:a.7} X(A_{\mathbf{q}}) 
= 2\,B_0\,A_{\mathbf{q}} -
n_3\, a_{3,2}\,A_{\mathbf{n}}\,B_2 - a_{2,1}\,A_{\mathbf{p}}.
\end{equation}
Since $A_\mathbf{n}$ and $A_\mathbf{m}$ are
constants, and $r=3$, we conclude that
$n_3\,a_{3,1}\,A_{\mathbf{n}} + 2\,a_{2,1}\,A_{\mathbf{m}} = 0$
and $A_{\mathbf{p}} = 0$. Otherwise, we should have $r = 2$ following the
similar discussion in (2). Let $A_{\mathbf{p}} =
0$ and $A_\mathbf{n} = 1$ in \eqref{eq:a.7}, we have
$$X(-\frac{A_{\mathbf{q}}}{n_3\,a_{3,2}}) = 2\,B_0\,(-\frac{A_\mathbf{q}}{n_3\,a_{3,2}}) + B_2.$$
From Lemma \ref{lem:a.15}, there exists a first integral $\omega$ of
$X$ such that
$$\dfrac{2\delta_2\omega\cdot\delta_2^3\omega-3\,(\delta_2^2\omega)^2}{(\delta_2\omega)^2} = -\frac{A_\mathbf{q}}{n_3\,a_{3,2}}\in K[\A_0(y)],$$
which implies (3) following the argument similar to the previous discussion.

Finally, we will show that it is impossible to have $3< r <\infty$. If on the contrary, we have $3 < r < \infty$, from Lemma
\ref{lem:a.11}, we can write $\mathbf{n} =
(n_1,0,n_3,\cdots,n_r)$. Let
\begin{eqnarray*}
\mathbf{m} &=& d_1^r(\mathbf{n}) =
(n_1,0,n_3,\cdots,n_{r-1}+1,n_r-1),\\
\mathbf{p} &=& {b_1^2}(\mathbf{m}) =
(n_1-1,1,n_3,\cdots,n_{r-1}+1,n_r-1),\\
\mathbf{q} &=& d_1^{r-1}(\mathbf{p}) =
(n_1-1,1,n_3,\cdots,n_{r-2}+1, n_{r-1},n_r-1),
\end{eqnarray*}
then $C(\mathbf{p}) = C(\mathbf{n})$ and therefore $\#(\mathbf{p})
= 0$ according to Lemma \ref{lem:a.10}. Let $\mathcal{J}_A$ to be the minimal closed
subsystem of $\mathcal{I}_A$ that contains $\mathbf{q}$, then $\mathbf{p}\in
\mathcal{J}_A$ and has the highest degree. However, $p_2 = 1\not=0$, which is
contradict to Lemma \ref{lem:a.11}. Thus, we concluded that either
$r = \infty$ or $r\leq 3$.
\end{proof}

\begin{proof}[Proof of Theorem \ref{th:8}]
The first part of Theorem \ref{th:8} is concluded from Lemma \ref{lem:a.7}
and Lemma \ref{le:2}.

It is easy to verify that when $r\leq 2$, the group $G$ is solvable, and when $r = 3$, $G$ is unsolvable.  Here we omit the detail calculations. In fact, we can consider $G$ as a multi-parameter group of transformation on $z\in \Comp$, and therefore the arguments for the corresponding multi-parameter Lie group are applicable. It is known that the transformation Lie groups with 1 or 2 parameters are solvable, while the three-parameter group is unsolvable\cite[pp. 85-86]{Blu:89}.  When $r=\infty$, $G$ contains an unsolvable subgroup that is defined by \eqref{eq:14}, and therefore is unsolvable.  
\end{proof}

\bibliographystyle{amsplain}

\subsection*{Acknowledgements}
The author is highly grateful to Keying Guan who leads 
him to the field of differential Galois theory. The author wishes
to express his hearty thanks to Claude Mitschi,
Universit\'{e} Louis Pasteur, Michael F. Singer, North
Carolina State University, and Colin Christopher,
University of Plymouth, for their great interested in reading the
early versions of the manuscript, and for their valuable advices
and helpful discussions.

\end{document}